\documentclass[12pt,a4paper]{article}

\usepackage{amssymb}

\usepackage{mathrsfs}
\usepackage{amsmath}

\usepackage{pstricks}

\newtheorem{theorem}{Theorem}[subsection]
\newtheorem{definition}[theorem]{Definition}
\newtheorem{lemma}[theorem]{Lemma}
\newtheorem{corollary}[theorem]{Corollary}

\newtheorem{example}[theorem]{Example}

\begin{document}
\title{Some new results on Gr\"{o}bner-Shirshov bases\footnote{Supported by the
NNSF of China (Nos. 10771077, 10911120389).}}
\author{
L. A. Bokut\footnote {Supported by RFBR 09-01-00157,
LSS--3669.2010.1 and SB RAS Integration grant No. 2009.97 (Russia)
and Federal Target Grant ¡°Scientific and educational personnel of
innovation Russia¡± for 2009-2013
(government contract No. 02.740.11.5191).} \\
{\small \ School of Mathematical Sciences, South China Normal
University}\\
{\small Guangzhou 510631, P. R. China}\\
{\small Sobolev Institute of Mathematics, Russian Academy of
Sciences}\\
{\small Siberian Branch, Novosibirsk 630090, Russia}\\
{\small bokut@math.nsc.ru}\\
\\
Yuqun Chen\footnote {Corresponding author.}
\\
{\small \ School of Mathematical Sciences, South China Normal
University}\\
{\small Guangzhou 510631, P. R. China}\\
{\small yqchen@scnu.edu.cn}\\
\\
K. P. Shum
\\
{\small \ Department of Mathematics, The
University of Hong Kong}\\
{\small Hong Kong, P. R. China}\\
{\small kpshum@maths.hku.hk}\\
}
\date{}
\maketitle

\begin{abstract}
In this survey article, we report some new results of
Gr\"{o}bner-Shirshov bases, including  new Composition-Diamond
lemmas and some  applications of some known Composition-Diamond
lemmas.
\end{abstract}

\noindent \textbf{Key words: }{Gr\"{o}bner-Shirshov basis;
Composition-Diamond lemma; normal form; Lie algebra; right-symmetry
algebra; Rota-Baxter algebra; $L$-algebra; conformal algebra;
metabelian Lie algebra, dendriform algebra; category; group; braid
group; inverse semigroup; plactic monoid; quantum group; operad.}

\noindent \textbf{AMS 2000 Subject Classification}: 16S15, 13P10,
13N10, 16-xx, 17-xx, 17B01, 17B66,  18Axx, 20F36, 20M18, 17B37,
18D50

\section{Introduction}

We continue our previous survey  \cite{survey08} on papers given
mostly by participants of an Algebra Seminar at South China Normal
University, Guangzhou.

This survey consists of  two blocks:
\begin{enumerate}
\item[(I)] \textit{New Composition-Diamond (CD-) lemmas.}
\item[(II)] \textit{Applications of known CD-lemmas.}
\end{enumerate}

In \cite{Sh62b}, A. I. Shirshov established the theory of
one-relator Lie algebras $Lie(X|s=0)$ and provided the algorithmic
decidality of the word problem for any one-relator Lie algebra. In
order to proceed his ideas, he first created the so called
Gr\"{o}bner-Shirshov bases theory for Lie algebras $Lie(X|S)$ which
are presented by generators and defining relations. The main notion
of Shirshov's theory was a notion of composition $(f,g)_w$ of two
Lie polynomials, $f,g\in Lie(X)$ relative to some associative word
$w$. The following lemma was proved by Shirshov \cite{Sh62b}.

\noindent{\bf Shirshov Lemma.} \emph{Let $Lie(X)\subset k\langle
X\rangle$ be a free Lie algebra over a field $k$ which is regarded
as an algebra of the Lie polynomials in the free algebra $k\langle
X\rangle$, and let $S$ be a Gr\"{o}bner-Shirshov basis in $Lie(X)$.
Then $f\in Id(S)\Rightarrow\bar{f}=u\bar{s}v$, where $s\in S^c, \
u,v\in X^*, \ \bar{f},\bar{s}$ are the leading associative words of
the corresponding Lie polynomials $f,s$ with respect to the deg-lex
order on $X^*$.}

Nowadays, Shirshov's lemma is named the ``Composition-Diamond lemma"
for Lie and associative algebras.

We explain  what it means of ``CD-lemma"  for a class (variety or
category) $\cal{M}$ of linear $\Omega$-algebras over a field $k$
(here $\Omega$ is a set of linear operations on $\cal{M}$) with free
objects.

\textbf{$\cal{M}$-CD-lemma} Let $\cal{M}$ be a class of (in general,
non-associative) $\Omega$-algebras, ${\cal M}(X)$ a free
$\Omega$-algebra in $\cal{M}$ generated by $X$ with a linear basis
consisting of ``normal (non-associative $\Omega$-) words" $[u]$, $S
\subset \ {\cal M}(X)$ a monic subset, $<$ a monomial order on
normal words and $Id(S)$ the ideal of ${\cal M}(X)$ generated by
$S$. Let $S$ be a Gr\"{o}bner-Shirshov basis (this means that any
composition of elements in $S$ is trivial). Then

\begin{enumerate}
\item[(a)] If $f\in Id(S)$, then $[\bar{f}]=[a\bar{s}b]$, where $[\bar{f}]$
is the leading word of $f$ and $[asb]$ is a normal $S$-word.

\item[(b)] $Irr(S)$ =$\{[u]| [u] \neq[a\bar{s}b], s\in S,
[asb] \ is \ a \ normal \ S-word\}$ is a linear basis of the
quotient algebra ${\cal M}(X|S)={\cal M}(X)/Id(S)$.
\end{enumerate}

In many cases, each of conditions (a) and (b) is equivalent to the
condition that $S$ is a Gr\"{o}bner-Shirshov basis in ${\cal M}(X)$.
But in some of our ``new CD-lemmas", this is not the case.

How to establish a CD-lemma for the free algebra ${\cal M}(X)$?
Following the idea of Shirshov  (see, for example,
\cite{BC07,BoChLi-Nankai,Sh62b}), one needs

1) to define appropriate linear basis (normal words) $N$ of ${\cal
M}(X)$;

2) to find a monomial order on $N$;

3) to find normal $S$-words;

4) to define compositions of elements in $S$ (they may be
compositions of intersection, inclusion and left (right)
multiplication, or may be else);

5) to prove  two key lemmas:

Key Lemma 1. Let $S$ be a Gr\"{o}bner-Shirshov basis (any
composition of polynomials from $S$ is trivial). Then  any element
of $Id(S)$  is a linear combination of normal $S$-words.

Key Lemma 2. Let $S$ be a Gr\"{o}bner-Shirshov basis, $[a_1s_1b_1]$
and $w=[a_2s_2b_2]$ normal $S$-words, $s_1, s_2\in S$. If
$\overline{[a_1s_1b_1]} = \overline{[a_2s_2b_2]}$, then $[a_1s_1b_1]
- [a_2s_2b_2]$ is trivial modulo $(S,w)$.

\ \

Let the normal words $N$ of the free algebra ${\cal M}(X)$ be a
well-ordered set and $0\neq f\in {\cal M}(X)$. Denote by $\bar f$
the leading word of $f$. $f$ is called monic if the coefficient of
$\bar f$ is 1.

A well ordering on $N$ is monomial if for any $u,v,w\in N$,
$$
u>v\Longrightarrow \overline{w|_{u}}>\overline{w|_{v}},
$$
where $w|_u=w|_{x\mapsto u}$ and $w|_v=w|_{x\mapsto v}$.

For example, let $X$ be a well-ordered set and $X^*$ the free monoid
generated by $X$. We define the deg-lex order on $X^*$: to compare
two words first by length and then lexicographically. Then, such an
order is monomial on  $X^*$.

Let $S\subset {\cal M}(X),\ s\in S,\ u\in N$. Then, roughly
speaking,  $u|_{s}=u|_{x_i\mapsto s}$, where $x_{i}$ is the
individuality occurrence of the letter $x_{i}\in X$ in $u$, is
called an $S$-word (or $s$-word). More precisely, let $\star\not\in
X$ be a new letter and $N(\star)$ the normal words of the free
algebra $ {\cal M}(X\cup\{\star\})$. Then by a $\star$-word we mean
any expression in $N(\star)$ with only one occurrence of $\star$.
Now, an $s$-word means $u|_{s}=u|_{\star\mapsto s}$, where $u$ is a
$\star$-word. Therefore,  the ideal $Id(S)$ of ${\cal M}(X)$ is the
set of linear combination of $S$-words.

An $S$-word $u|_s$ is normal if $\overline{u|_s}=u|_{\overline{s}}$.

\begin{definition}
Given a monic subset $S\subset {\cal M}(X)$ and $w\in N $, an
intersection (inclusion) composition  $h$ is called trivial modulo
$(S,w)$ if $h$ can be presented as a linear combination of normal
$S$-words with leading words less than  $w$; a left (right)
multiplication composition (or another kinds composition) $h$ is
called trivial modulo $(S)$ if $h$ can be presented as a linear
combination of normal $S$-words with leading words less than or
equal to $\bar h$.

The set $S$ is a Gr\"{o}bner-Shirshov basis in ${\cal M}(X)$ if all
the possible compositions of elements in $S$ are trivial modulo $S$
and corresponding $w$.
\end{definition}

If a subset $S$ of ${\cal M}(X)$ is not a Gr\"{o}bner-Shirshov basis
then one can add all nontrivial compositions of polynomials of $S$
to $S$. Continuing this process repeatedly, we finally obtain a
Gr\"{o}bner-Shirshov basis $S^{C}$ that contains $S$. Such a process
is called Shirshov algorithm. $S^{C}$ is called Gr\"{o}bner-Shirshov
complement of $S$.

\ \

We establish new CD-lemmas for the following classes of algebras:

\ \
\begin{enumerate}
\item[2.1]\
Lie algebras over commutative algebras (L. A. Bokut, Yuqun Chen,
Yongshan Chen \cite{BoChChen-Liecomm}).

\item[2.2]\ Associative algebras with multiple operations (L. A. Bokut,
Yuqun Chen, Jianjun Qiu \cite{BoChQiu-CD-Omega}).

\item[2.3]\ Rota-Baxter
algebras (L. A. Bokut, Yuqun Chen, Xueming Deng
\cite{BoChDeng-RotaB}).

\item[2.4]\  Right-symmetric (or pre-Lie) algebras (L. A. Bokut,
Yuqun Chen, Yu Li \cite{BoChLi-GSB-rightsym}).

\item[2.5]\ $L$-algebras
(Yuqun Chen, Jiapeng Huang \cite{ChHuang-L-alg}).

\item[2.6]\ Differential
algebras (Yuqun Chen, Yongshan Chen, Yu Li \cite{ChChLi-GSB-diff}).

\item[2.7]\ Non-associative algebras  over commutative algebras (Yuqun
Chen, Jing Li, Mingjun Zeng \cite{ChLiZeng-Cd-nonasso}).

\item[2.8]\
$n$-conformal algebras (L. A. Bokut, Yuqun Chen, Guangliang Zhang
\cite{BoChZhang-n-conf}).

\item[2.9]\ $\lambda$-differential  associative
algebras with multiple operators (Jianjun Qiu, Yuqun Chen
\cite{ChQiu-Cd-diff}).

\item[2.10]\ Categories (L. A. Bokut, Yuqun Chen, Yu
Li \cite{BoChLi-Nankai}).

\item[2.11]\  Metabelian Lie algebras  (Yuqun
Chen, Yongshan Chen \cite{ChCh-Matabel}).

\item[2.12]\  $S$-act algebras  (Xia Zhang \cite{Xia}).

\ \

A new CD-lemma is given in:

\item[2.13]\ Operads  (V.
Dotsenko, A. Khoroshkin \cite{DK10}, see also V. Dotsenko, M. V.
Johansson \cite{DJo}).
\end{enumerate}

\ \

Some new applications of known CD-lemmas are given in:

\ \

\begin{enumerate}
\item[3.1]\ Artin-Markov normal form for braid
groups in Artin-Burau generators (L. A. Bokut,  V. V. Chainikov, K.
P. Shum \cite{BoChaSh07} and Yuqun Chen, Qiuhui Mo
\cite{ChMo-GSB-Braid-Artin}).

\item[3.2]\ Braid
groups in Artin-Garside generators (L. A. Bokut \cite{Bo08}).

\item[3.3]\ Braid
groups in Birman-Ko-Lee generators (L. A. Bokut \cite{Bo09}).

\item[3.4]\ Braid
groups in Adyan-Thurston generators (Yuqun Chen, Chanyan Zhong
\cite{ChZhong-Braid}).

\item[3.5]\
GSB and normal forms for free inverse semigroups (L. A. Bokut, Yuqun
Chen, Xiangui Zhao \cite{BoChZhao-inverse-sg}).

\item[3.6]\ Embeddings of algebras
(L. A. Bokut, Yuqun Chen, Qiuhui Mo \cite{BoChMo-embed}).

\item[3.7]\
Word problem for Novikov's and Boone's group (Yuqun Chen, Wenshu
Chen, Runai Luo \cite{ChChLuo}).

\item[3.8]\  PBW basis of $U_q(A_N)$ (L. A. Bokut, P. Malcolmson \cite{BoMa} and Yuqun Chen,
Hongshan Shao, K. P. Shum \cite{ChShaoShum}).

\item[3.9]\ GSB for free dendriform algebras (Yuqun Chen, Bin Wang
\cite{ChWang-GSB-Den}).

\item[3.10]\
Anti-commutative GSB of a free Lie algebra relative to
Lyndon-Shirshov words (L. A. Bokut, Yuqun Chen, Yu Li
\cite{BoChLi10}).

\item[3.11]\ Free
partially commutative groups, associative algebras and Lie algebras
(Yuqun Chen, Qiuhui Mo \cite{ChMo-Partialcomm}).

\item[3.12]\ Plactic monoids in row generators (Yuqun Chen, Jing Li
\cite{ChLiJing-plactic}).

\ \

Some applications of CD-lemmas for associative and Lie algebras are
given in:

\item[3.13]\ Plactic algebras in standard generators (Lukasz Kubat, Jan Okni\'{n}skii \cite{Okninskii}).

\item[3.14]\ Filtrations and distortion in
infinite-dimensional algebras (Y. Bahturin, A. Olshanskii
\cite{Olshanskii}).

\item[3.15]\ Sufficiency conditions for Bokut' normal forms (K. Kalorkoti
\cite{Kalorkoti10}).

\item[3.16]\ Quantum groups of type $G_2$ and $D_4$ (Yanhua Ren,
Abdukadir Obul \cite{Abdu1}, Gulshadam Yunus, Abdukadir Obul
\cite{Abdu2}).

\item[3.17]\ An embedding of recursively presented Lie algebras (E. Chibrikov
\cite{Chibrikov}).

\item[3.18]\ GSB for some monoids (Canan Kocapinar, Firat Ates, A. Sinan Cevik \cite{Turkey10}).
\end{enumerate}

\ \

Notations

\ \

CD-lemma: Composition-Diamond lemma.

GSB: Gr\"{o}bner-Shirshov basis.

$X^*$: the free monoid generated by the set $X$.

$X^{\ast \ast }$: the set of all non-associative words $(u)$ in $X$.

$k$: a field.

$K$:  a commutative ring with unit.

$Id(S)$: the ideal generated by the set $S$.

$k\langle X\rangle$: the free associative algebra over $k$ generated
by $X$.

$k[ X]$: the polynomial algebra over $k$ generated by $X$.

$Lie(X)$: the free Lie algebra over $k$ generated by $X$.

$Lie_K(X)$: the free Lie algebra over a commutative ring $K$
generated by $X$.

\section{New CD-lemmas}

\subsection{Lie algebras over a commutative algebra}

A. A. Mikhalev, A. A. Zolotykh \cite{MZ} prove the CD-lemma for
$k[Y]\otimes k\langle X\rangle$, a tensor product of a free algebra
and a polynomial algebra.  L. A. Bokut, Yuqun Chen, Yongshan Chen
\cite{BCC08} prove the CD-lemma for $k\langle Y\rangle\otimes
k\langle X\rangle$, a tensor product of two free algebras. Yuqun
Chen, Jing Li, Mingjun Zeng \cite{ChLiZeng-Cd-nonasso} (see \S2.7)
prove the CD-lemma for $k[Y]\otimes k(X)$, a tensor product of a
non-associative algebra and a polynomial algebra.

In this subsection, we introduce the CD-lemma for $k[Y]\otimes Lie(
X)$, Lie algebra $Lie(X)$ over a polynomial algebra $k[Y]$, which is
established in \cite{BoChChen-Liecomm}. It provides a
Gr\"{o}bner-Shirshov bases theory for Lie algebras over a
commutative algebra.

Let $K$ be a commutative associative ${k}$-algebra with unit and
$\mathcal{L}$ a Lie $K$-algebra.  Then,  $\mathcal{L}$ can be
presented as $K$-algebra by generators $X$ and some defining
relations $S$,
$$
                              \mathcal{L}=Lie_K(X|S)=Lie_K(X)/Id(S).
$$
In order to define a Gr\"{o}bner-Shirshov basis for $\mathcal{L}$,
we first present $K$ in a form

$$
                           K={k}[Y|R]={k}[Y]/Id(R),
$$
where  $R\subset{k}[Y]$. Then the Lie $K$-algebra $\mathcal{L}$ has
the following presentation as a ${k}[Y]$-algebra
$$
          \mathcal{L}= Lie_{{k}[Y]}(X| S, Rx,\ x\in X).
$$

\ \

Let the set $X$ be  well-ordered, and let $<$ and $\prec_X$ be the
lex order and the deg-lex order
 on $X^*$ respectively, for example, $ab<a,\ a\prec_X ab,\ a,b\in X$.

A word $w\in X^*\setminus\{1\}$ is an associative Lyndon-Shirshov
 word (ALSW for short) if
 $$
 (\forall u,v\in X^*, u,v\neq1) \ w=uv\Rightarrow w>vu.
 $$

 A nonassociative word $(u)\in X^{**}$  is a non-associative
Lyndon-Shirshov word (NLSW for short), denoted by $[u]$,  if

(i) $u$ is an ALSW;

(ii) if $[u]=[(u_{1})(u_{2})]$ then both $(u_{1})$ and $(u_{2})$ are
NLSW's (from (i) it then follows that $u_{1}>u_{2}$);

(iii) if $[u]=[[[u_{11}][u_{12}]][u_{2}]]$ then $u_{12} \leq u_{2}$.

We denote the set of all NLSW's (ALSW's) on $X$ by $NLSW(X)$\
($ALSW(X)$).

Let $Y=\{y_j|j\in J\}$ be a well-ordered set and
$[Y]=\{y_{j_1}y_{j_2}\cdots y_{j_l}|y_{j_1}\leq y_{j_2}\leq\cdots
\leq y_{j_l}, l \geq0 \}$ the free commutative monoid generated by
$Y$. Then $[Y]$ is a ${k}$-linear basis of the polynomial algebra
${k}[Y]$.

Let $Lie_{{k}[Y]}(X)$ be the ``double" free Lie algebra, i.e., the
free Lie algebra over the polynomial algebra ${k}[Y]$ with
generating set $X$.

From now on we regard $Lie_{{k}[Y]}(X)\cong {k}[Y]\otimes Lie_{{
k}}(X)$ as the Lie subalgebra of ${k}[Y] \langle X \rangle\cong{\bf
k}[Y]\otimes{k} \langle X \rangle$ the free associative algebra over
polynomial algebra ${k}[Y]$, which is generated by $X$ under the Lie
bracket $[u,v]=uv-vu$.

For an $ALSW$ $t$ on $X$, there is a unique bracketing, denoted by
$[t]$, such that $[t]$ is $NLSW$:
$$
[t]=t\in X, \ [t]=[[v][w]],
$$
where $t=vw$ and $w$ is the longest $ALSW$  proper end of $t$ (then
$v$ is also an ALSW).

Let
$$
T_A=\{u=u^Yu^X|u^Y\in[Y], \ u^X  \in ALSW(X)\}
$$
and
$$
T_N=\{[u]=u^Y[u^X]|u^Y\in[Y], \ [u^X]\in NLSW (X)\}.
$$

Let ${k}T_N$ be the linear space spanned by $T_N$ over ${k}$. For
any $[u],[v]\in T_N$, define
$$
[u][v]=\sum\alpha_iu^Yv^Y[w_i^X]
$$
where $\alpha_i\in {k},\ [w_i^X]$'s are NLSW's and
$[u^X][v^X]=\sum\alpha_i[w_i^X]$ in $Lie_{{k}}(X)$.

Then ${k}[Y]\otimes Lie_{{k}}(X)\cong {k}T_N$ as ${k}$-algebra and
$T_N$ is a ${k}$-basis of ${k}[Y]\otimes Lie_{{k}}(X)$.

We define the deg-lex order $\succ$ on
$$[Y]X^*=\{u^Yu^X|u^Y\in
[Y],u^X\in X^*\}
$$
by the following: for any $u,v\in [Y]X^*$,
\begin{eqnarray*}
&&u\succ v  \ \mbox{ if } \ (u^X\succ_Xv^X)  \  \mbox{ or }\
(u^X=v^X \  \mbox{ and }\ u^Y\succ_Yv^Y ),
\end{eqnarray*}
where $\succ_Y$ and $\succ_X$ are the deg-lex order on $[Y]$ and
$X^*$ respectively.

\begin{lemma}\label{l2}(Shirshov \cite{Shir1, Shir3}) Suppose that
$w=aub$ where $w,u\in T_A$ and $a,b\in X^*$. Then
$
[w]=[a[uc]d],
$
where $[uc]\in T_N$  and $b=cd$.

Represent $c$ in a form $c=c_{1}c_{2} \ldots c_{k},$ where $c_{1},
\ldots ,c_{n}\in ALSW(X)$ and $c_{1} \leq c_{2} \leq \ldots \leq
c_{n}$. Denote by
$$
[w]_u=[a[\cdots[[[u][c_{1}]][c_{2}] ]\ldots [c_{n}]]d].
$$
Then, $ \overline{[w]}_{u}=w. $
\end{lemma}

\begin{definition} Let $S\subset Lie_{{k}[Y]}(X)$ be a ${k}$-monic subset, $a,b\in X^*$ and
$s\in S$. If $a\bar{s}b\in T_A$,  then by Lemma \ref{l2} we have the
special bracketing $[a\bar{s}b]_{\bar{s}}$ of $a\bar{s}b$ relative
to  $\bar{s}$.  We define
$[asb]_{\bar{s}}=[a\bar{s}b]_{\bar{s}}|_{[\bar{s}]\mapsto{s}}$ to be
a normal $s$-word (or normal $S$-word).
\end{definition}

It is proved that each element in $Id(S)$ can be expressed as a
linear combinations of normal $S$-words (see
\cite{BoChChen-Liecomm}). The proof is not easy.

There are four kinds of compositions.

\begin{definition}\label{d1}
Let $f,g$ be  two ${k}$-monic polynomials of $Lie_{{k}[Y]}(X)$.
Denote the least common multiple of $\bar{f}^Y$ and $\bar{g}^Y$ in
$[Y]$ by $L=lcm(\bar{f}^Y,\bar{g}^Y)$.

If $\bar{g}^X$ is a subword of $ \bar{f}^X$, i.e.,
$\bar{f}^X=a\bar{g}^Xb$ for some $a,b\in X^*$, then the polynomial
$$
C_1\langle
f,g\rangle_w=\frac{L}{\bar{f}^Y}f-\frac{L}{\bar{g}^Y}[agb]_{\bar{g}}
$$
is called the { inclusion composition} of $f$ and $g$ with respect
to $w$, where $w=L\bar{f}^X=La\bar{g}^Xb$.

If a proper prefix of $\bar{g}^X$ is a proper suffix of $\bar{f}^X$,
i.e., $\bar{f}^X=aa_0$, $\bar{g}^X=a_0b$, $a,b,a_0\neq1$, then the
polynomial
$$
C_2\langle
f,g\rangle_w=\frac{L}{\bar{f}^Y}[fb]_{\bar{f}}-\frac{L}{\bar{g}^Y}[ag]_{\bar{g}}
$$
is called the  {intersection composition} of $f$ and $g$ with
respect to $w$, where $w=L\bar{f}^Xb=La\bar{g}^X$.

 If the greatest common divisor of $\bar{f}^Y$ and $\bar{g}^Y$ in
$[Y]$ is non-empty, then for any $a,b,c\in X^*$ such that
$w=La\bar{f}^Xb\bar{g}^Xc\in T_A$, the polynomial
$$
C_3\langle
f,g\rangle_w=\frac{L}{\bar{f}^Y}[afb\bar{g}^Xc]_{\bar{f}}-
\frac{L}{\bar{g}^Y}[a\bar{f}^Xbgc]_{\bar{g}}
$$
is called the  {external composition} of $f$ and $g$ with respect to
$w$.

 If  $\bar{f}^Y\neq1$,
then for any normal $f$-word $[afb]_{\bar{f}}, \ a,b\in X^*$, the
polynomial
$$
C_4\langle f\rangle_w=[a\bar{f}^Xb][afb]_{\bar{f}}
$$
is called the  {multiplication composition} of $f$ with respect to
$w$, where $w=a\bar{f}^Xba\bar{f}b$.
\end{definition}

\begin{theorem}{\bf (\cite{BoChChen-Liecomm},
CD-lemma for Lie algebras over a commutative algebra) }
\label{cdL2.1.4} Let $S\subset{Lie_{{k}[Y]}(X)}$ be nonempty set of
${k}$-monic polynomials and $Id(S)$ be the ${k}[Y]$-ideal of
$Lie_{{k}[Y]}(X)$ generated by $S$. Then the following statements
are equivalent.
\begin{enumerate}
\item[(i)] $S$ is a Gr\"{o}bner-Shirshov basis in
$Lie_{{k}[Y]}(X)$.
\item[(ii)] $f\in{Id(S)}\Rightarrow{\bar{f}=\beta a\bar{s}b\in T_A}$ for
some $s\in{S},  \ \beta\in[Y]$ and $a,b\in{X^*}$.
\item[(iii)]$Irr(S)=\{[u] \ | \ [u]\in T_N, \ u\neq{\beta a\bar{s}b},  \mbox{ for any }
s\in{S},\ \beta\in[Y], \ a,b\in{X^*}\}$ is a $k$-basis for
$Lie_{{k}[Y]}(X|S)=Lie_{{k}[Y]}(X)/Id(S)$.
\end{enumerate}
\end{theorem}

Now, we give some applications of Theorem \ref{cdL2.1.4}.

A Lie algebra is called special if it can be embedded into its
universal enveloping associative algebra.

Only a few of non-special Lie algebras were known.

\begin{example}(Shirshov \cite{Shir53,Shir3})
Let the field ${k}=GF(2)$ and $K={k}[Y|R]$, where
$$Y=\{y_i,i=0,1,2,3\}, \ R=\{y_0y_i=y_i \ (i=0,1,2,3), \ y_iy_j=0 \ (i,j\neq0)\}.$$
Let $\mathcal{L}=Lie_K(X|S_1, S_2)$, where $X=\{x_i, 1\leq i\leq
13\},$ $S_1$ consist of the following relations
\begin{eqnarray*}
&&[x_2,x_1]=x_{11}, \ [x_3,x_1]=x_{13}, \ [x_3,x_2]=x_{12}, \\
&&[x_5,x_3]=[x_6,x_2]=[x_8,x_1]=x_{10}, \\
&& [x_i,x_j]=0 \ \ \ (\mbox{for any other} \ i>j),
\end{eqnarray*}
and $S_2$ consist of the following relations
\begin{eqnarray*}
&& y_0x_i=x_i  \ (i=1,2,\ldots,13),\\
&& x_4=y_1x_1, \ x_5=y_2x_1, \ x_5=y_1x_2, \ x_6=y_3x_1, \
x_6=y_1x_3, \\
&& x_7=y_2x_2, \ x_8 =y_3x_2, \ x_8 =y_2x_3,  \ x_9=y_3x_3, \\
&& y_3x_{11}=x_{10}, \ y_1x_{12}=x_{10}, \ y_2x_{13}=x_{10},\\
&& y_1x_k=0 \ (k=4,5,\ldots,11,13), \ y_2x_t=0 \ (t=4,5,\ldots,12),
\ y_3x_l=0 \ (l=4,5,\ldots,10,12,13).
\end{eqnarray*}
Then $S=S_1\cup S_2\cup RX\cup\{ y_1x_2=y_2x_1, \ y_1x_3=y_3x_1, \
y_2x_3=y_3x_2\}$ is a GSB in $Lie_{{\bf k}[Y]}(X)$. By Theorem
\ref{cdL2.1.4}, $x_{10}\neq0$ in $\mathcal{L}$. But from the
CD-lemma for $k[Y]\otimes k\langle X\rangle$ in  A. A. Mikhalev, A.
A. Zolotykh \cite{MZ}, it easily follows that $x_{10}=0$ in the
universal enveloping associative algebra $U_K(\mathcal{L})=K\langle
X|S_1^{(-)},S_2\rangle$. Hence, $\mathcal{L}$ is non-special as a
$K$-algebra.
\end{example}

\begin{example} (Cartier  \cite{Cartier}) Let ${k}=GF(2)$,
$K={k}[y_1,y_2,y_3|y_i^2=0,\ i=1,2,3]$ and
$\mathcal{L}=Lie_{K}(X|S)$, where $X=\{x_{ij},1\leq i\leq j\leq3\}$
and
$$
S=\{[x_{ii},x_{jj}]=x_{ji} \ (i>j), [x_{ij},x_{kl}]=0 \
(\mbox{others}), \ y_3x_{33}=y_2x_{22}+y_1x_{11}\}.
$$
Then $S'=S\cup \{y_i^2x_{kl}=0\ (\forall i,k,l)\}\cup S_1$ is a GSB
in $Lie_{{k}[Y]}(X)$, where $S_1$ consists of the following
relations
\begin{eqnarray*}
&&y_3x_{23}=y_1x_{12}, \  y_3x_{13}=y_2x_{12}, \ y_2x_{23}=y_1x_{13}, \ y_3y_2x_{22}=y_3y_1x_{11}, \\
&&y_3y_1x_{12}=0, \ y_3y_2x_{12}=0,  \ y_3y_2y_1x_{11}=0,   \
y_2y_1x_{13}=0.
\end{eqnarray*}

The universal enveloping algebra of $\mathcal{L}$ has a
presentation:
$$
U_K(\mathcal{L})=K\langle X|S^{(-)}\rangle\cong{\bf k}[Y]\langle
X|S^{(-)},y_i^2x_{kl}=0\ (\forall i,k,l)\rangle.
$$

In $U_K(\mathcal{L})$, we have
\begin{eqnarray*}
0=y_3^2x_{33}^2=(y_2x_{22}+y_1x_{11})^2=y_2^2x_{22}^2+y_1^2x_{11}^2+y_2y_1[x_{22},x_{11}]
 = y_2y_1x_{12}.
\end{eqnarray*}
On the other hand, since $y_2y_1x_{12}\in Irr(S')$,
$y_2y_1x_{12}\neq0$ in $\mathcal{L}$ by Theorem \ref{cdL2.1.4}.
Thus, $\mathcal{L}$ is not special as a $K$-algebra.
\end{example}

\noindent{\bf Conjecture} (Cohn \cite{Conh63}) Let $K={\bf
k}[y_1,y_2,y_3|y_i^p=0, i=1,2,3]$ be an algebra of truncated
polynomials over a field $k$ of characteristic $p>0$. Let
$$
\mathcal{L}_p=Lie_{K}(x_1,x_2,x_3 \ | \ y_3x_3=y_2x_2+y_1x_1).
$$
Then $\mathcal{L}_p$ is not special. $\mathcal{L}_p$ is called
Cohn's Lie  algebras.

Let $Y=\{y_1,y_2,y_3\}$, $X=\{x_1,x_2,x_3\}$ and
$S=\{y_3x_3=y_2x_2+y_1x_1, \  y_i^px_j=0,\ 1\leq i,j\leq3\}$. Then
$\mathcal{L}_p\cong Lie_{{\bf k}[Y]}(X  | S)$ and
 $U_K(\mathcal{L}_p)\cong{\bf k}[Y]\langle X|S^{(-)}\rangle$.
 Suppose that $S^C$ is the Gr\"{o}bner-Shirshov complement
 of $S$ in $Lie_{{\bf k}[Y]}(X)$. Let $S_{_{X^{^p}}}\subset\mathcal{L}_p$
 be the set of all the elements of
$S^C$ whose $X$-degrees do not exceed $p$. It is clear that
$S_{_{X^{^p}}}$ is a finite set for any $p$. Although it is
difficult to find the Gr\"{o}bner-Shirshov complement
 of $S$ in $Lie_{{\bf k}[Y]}(X)$, it is possible to find the set $S_{_{X^{^p}}}$. By using this
idea and the Theorem \ref{cdL2.1.4}, we have the following theorem:

\begin{theorem}(\cite{BoChChen-Liecomm})\
Cohn's Lie  algebras $\mathcal{L}_2$, $\mathcal{L}_3$ and
$\mathcal{L}_5$ are non-special.
\end{theorem}

\begin{theorem}\label{t4.5}(\cite{BoChChen-Liecomm})\  For an arbitrary
commutative ${k}$-algebra $K={k}[Y|R]$, if $S$ is a
Gr\"{o}bner-Shirshov basis in $Lie_{{k}[Y]}(X)$ such that for any
$s\in S$, $s$ is ${k}[Y]$-monic, then $\mathcal{L}=Lie_{K}(X|S)$ is
special.
\end{theorem}

\begin{corollary}\label{co4.6}(\cite{BoChChen-Liecomm})\
Any Lie $K$-algebra $L_K=Lie_K(X|f)$ with one monic defining
relation $f=0$ is special.
\end{corollary}

\begin{theorem}(\cite{BoChChen-Liecomm})\ Suppose that $S$ is a
finite homogeneous subset of $Lie_{{k}}(X)$. Then the word problem
of $Lie_{K}(X|S)$ is solvable for any finitely generated commutative
${k}$-algebra $K$.
\end{theorem}

\begin{theorem}(\cite{BoChChen-Liecomm})\ Every finitely or
countably generated Lie $K$-algebra can be embedded into a two-generated Lie
$K$-algebra, where $K$ is an arbitrary commutative ${\bf
k}$-algebra.
\end{theorem}

\subsection{Associative algebras with multiple operations}

A $\Omega$-algebra $A$ is a $k$-space with the linear operator set
$\Omega$ on $A$.

V. Drensky and R. Holtkamp \cite{Dren} constructed Gr\"obner bases
theory for $\Omega$-algebras, where $\Omega$ consists of $n$-ary
operations, $n\geq2$.

In this subsection, we consider associative $\Omega$-algebras, where
$\Omega$ consists of $n$-ary operations, $n\geq1$.

An associative algebra with multiple linear operators is an
associative $K$-algebra $R$ with a set
 $\Omega$ of multilinear operators (operations).

Let $X$ be a set and
$$
\Omega=\bigcup_{n=1}^{\infty}\Omega_{n}
 $$
where $\Omega_{n}$ is the set of $n$-ary operations, for example,
ary $(\delta)=n$ if $\delta\in \Omega_n$.

Define
\begin{eqnarray*}
&&\mathfrak{S}_{0}=S(X_0),\  X_0=X,\\
&&\mathfrak{S}_{1}=S(X_1), \  X_1=X\cup \Omega(\mathfrak{S}_{0}), \\
&&\cdots\cdots\\
&& \mathfrak{S}_{n}=S(X_n), \ X_n=X\cup \Omega(\mathfrak{S}_{n-1}),
\end{eqnarray*}
where $
\Omega(\mathfrak{S}_{j})=\bigcup\limits_{t=1}^{\infty}\{\delta(u_1,u_2,\dots,u_t)|\delta\in
\Omega_t, u_i\in \mathfrak{S}_{j} , \ i=1,2,\dots,t\}$ and $ S(X_j)$
is the free semigroup generated by $X_j,\ j=0,1,\dots$.

Let
$$
\mathfrak{S}(X)=\bigcup_{n\geq0}\mathfrak{S}_{n}.
$$
Then, it is easy to see that $\mathfrak{S}(X)$ is a semigroup such
that $ \Omega(\mathfrak{S}(X))\subseteq \mathfrak{S}(X). $

 For any $u\in \mathfrak{S}(X)$, $dep(u)=\mbox{min}\{n|u\in\mathfrak{S}_{n} \}$
 is called the depth of
 $u$.

Let   $K\langle X; \Omega\rangle$ be the $K$-algebra spanned by
$\mathfrak{S}(X)$.  Then, the element in $\mathfrak{S}(X)$ (resp.
$K\langle X; \Omega\rangle$) is called a $\Omega$-word (resp.
$\Omega$-polynomial). If $u\in X\cup \Omega(\mathfrak{S}(X))$, we
call $u$ a prime $\Omega$-word and define $bre(u)=1$ (the breadth of
$u$). If $u=u_1u_2\cdots u_n\in\mathfrak{S}(X)$, where $u_i$ is
prime $\Omega$-word for all $i$, then we define $bre(u)=n$.

Extend linearly each $\delta\in\Omega_n$,
$$\delta:\mathfrak{S}(X)^n\rightarrow \mathfrak{S}(X), \
(x_1,x_2,\cdots,x_n)\mapsto \delta(x_1,x_2,\cdots,x_n)
$$
to $K\langle X; \Omega\rangle$. Then, $K\langle X; \Omega\rangle$ is
a free associative algebra with multiple linear operators $\Omega$
on set $X$.

Assume that $ \mathfrak{S}(X)$ is equipped with a monomial order
$>$.

Let $f, g$ be two monic $\Omega$-polynomials. Then, there are two
kinds of $compositions$.
\begin{enumerate}
\item[(i)]If there exists a $\Omega$-word $w=\bar{f}a=b\bar{g}$ for some $a,b\in
\mathfrak{S}(X)$ such that $bre(w)< bre(\bar{f})+bre(\bar{g})$, then
we call $(f,g)_{w}=fa-bg$ the $intersection$ $composition$ of $f$
and $g$ with respect to $w$.
\item[(ii)] If there exists a $\Omega$-word $w=\bar{f}=u|_{\bar{g}}$ for some
$u \in \mathfrak{S}(X)$, then we call $(f,g)_{w}=f-u|_{g}$ the
$including$ $composition$ of  $f$ and $g$ with respect to $w$.
\end{enumerate}

It is noted that each $S$-word is normal.

\begin{theorem}\label{CD-RotaB}{{\bf(\cite{BoChQiu-CD-Omega},
CD-lemma for associative $\Omega$-algebras)}}\ \  Let $S$ be a set
of monic $\Omega$-polynomials in $K\langle X;\Omega\rangle$ and
 $>$ a monomial order on $\mathfrak{S}(X)$.  Then the following
statements are equivalent.
 \begin{enumerate}
\item[(i)] $S $ is a Gr\"{o}bner-Shirshov basis in $K\langle X;\Omega\rangle$.
\item[(ii)] $ f\in Id(S)\Rightarrow \bar{f}=u|_{\overline{s}}$
for some $u \in \mathfrak{S}(X)$ and $s\in S$.
\item[(iii)] $Irr(S) = \{ w\in \mathfrak{S}(X) |  w \neq
u|_{\overline{s}}
 \mbox{ for  any} \ u \in \mathfrak{S}(X) \ \mbox{and } s\in S\}$
is a $K$-basis of $K\langle X;\Omega|S\rangle=K\langle
X;\Omega\rangle/Id(S)$.
\end{enumerate}
\end{theorem}

Now, we give some applications of Theorem \ref{CD-RotaB}.

First of all,  we define an order on $\mathfrak{S}(X)$. Let $X$ and
$\Omega$ be well-ordered sets. We order $X^*$ by the deg-lex order.
For any $u\in \mathfrak{S}(X)$, $u$ can be uniquely expressed
without brackets as
$$
u=u_0\delta_{i_{_{1}}}\overrightarrow{x_{i_1}}u_1\cdots
\delta_{i_{_{t}}}\overrightarrow{x_{i_{t}}}u_t,
$$
where each $u_i\in X^*,\delta_{i_{_{k}}}\in \Omega_{i_{_{k}}}, \
\overrightarrow{x_{i_k}}=( x_{k_{1}},x_{k_{2}},\cdots,
x_{k_{i_{_{k}}}} )\in \mathfrak{S}(X)^{i_k}$. It is reminded that
for each $i_{{k}}, \ dep(u)>dep(\overrightarrow{x_{i_k}})$.

Denote by
$$
wt(u)=(t,\delta_{i_{_{1}}},\overrightarrow{x_{i_{1}}},
\cdots,\delta_{i_{_{t}}},\overrightarrow{x_{i_t}}, u_0, u_1, \cdots,
u_t ).
$$
Then, we order  $\mathfrak{S}(X)$ as follows: for any $u,v \in
\mathfrak{S}(X) $,
\begin{equation}\label{o1}
u>v\Longleftrightarrow wt(u)>wt(v)\ \mbox{ lexicographically}
\end{equation}
by induction on $dep(u)+dep(v)$.

It is clear that the order (\ref{o1}) is a monomial order on
$\mathfrak{S}(X)$.

A Rota-Baxter $K$-algebra of weight $\lambda$ (\cite{Bax60, Ro69})
is an associative algebra $R$ with  a $K$-linear operation $
P:R\rightarrow R$ satisfying the Rota-Baxter relation:
$$
 P(x)P(y) = P( P(x)y +xP(y)+\lambda xy), \  \forall x,y \in  R.
$$

Thus, any Rota-Baxter algebra is a special case of associative
algebra with multiple operators when  $\Omega=\{P\}$.

Now, let  $\Omega=\{P\}$ and $\mathfrak{S}(X)$ be as before. Let
$K\langle X;P\rangle$ be the free associative algebra with one
operator $\Omega=\{P\}$ on a set $X$.

\begin{theorem}(\cite{BoChQiu-CD-Omega})\label{t2.2}\
With the order (\ref{o1})  on $\mathfrak{S}(X)$,
$$
S=\{P(x)P(y) - P( P(x)y) - P(xP(y))-\lambda P(xy) |\  x,y \in
\mathfrak{S}(X)\}
$$
is a Gr\"{o}bner-Shirshov basis in $K\langle X;P\rangle$.
\end{theorem}

By Theorems \ref{CD-RotaB} and \ref{t2.2}, we obtain a normal form
$Irr(S)$ of  the free Rota-Baxter algebra $K\langle X;P|S\rangle$
which is the same as in \cite{EG08a}.

\ \

A $\lambda$-differential
 algebra over $K$ (\cite{GK08})  is an  associative $K$-algebra  $R$
  together with a  $K$-linear operator $D:R\rightarrow R$ such that
$$
D(xy)=D(x)y+xD(y)+\lambda D(x)D(y),\ \forall x, y \in R.
$$

Any  $\lambda$-differential  algebra is also an   associative
algebra with one
 operator $\Omega=\{D\}$.

Let $X$ be well-ordered and $K\langle X;D\rangle$ the free
associative algebra with one operator  $\Omega=\{D\}$ defined as
before.

 For any $u\in
\mathfrak{S}(X)$, $u$ has a unique expression
$$
u=u_1u_2\cdots u_n,
$$
where each $u_i\in X\cup D(\mathfrak{S}(X))$. Denote by
$deg_{_{X}}(u)$ the number of $x\in X$ in $u$.   Let
$$
wt(u)=(deg_{_{X}}(u), u_1,u_2,\cdots, u_n).
$$
For any $u,v \in \mathfrak{S}(X) $, define
\begin{equation}\label{o2}
u>v\Longleftrightarrow wt(u)>wt(v)\ \mbox{ lexicographically},
\end{equation}
where for each  $t, \ u_t>v_t$ if one of the following holds:

(a) $u_t, v_t\in X$ and $u_t>v_t$;

(b)  $u_t=D(u_{t}^{'}), v_t\in X$;

(c)  $u_t=D(u_{t}^{'}),v_t=D(v_{t}^{'})$ and $u_{t}^{'}>v_{t}^{'}$.

\ \

It is easy to see the order (\ref{o2}) is a monomial order on
$\mathfrak{S}(X)$.

\begin{theorem}(\cite{BoChQiu-CD-Omega})\label{t2.2.3} With  the  order
(\ref{o2})  on $\mathfrak{S}(X)$,
$$
S=\{D(xy)-D(x)y-xD(y)-\lambda D(x)D(y) |\  x,y \in \mathfrak{S}(X)\}
$$
is a Gr\"{o}bner-Shirshov basis in $K\langle X;D\rangle$.
\end{theorem}

By Theorems \ref{CD-RotaB} and \ref{t2.2.3}, we obtain a normal form
$Irr(S)$ of  the free $\lambda$-differential  algebra $K\langle
X;D|S\rangle$ which is the same as in \cite{GK08}.

\ \

A differential Rota-Baxter algebra of weight $\lambda$
 (\cite{GK08}), called also $\lambda$-differential Rota-Baxter
 algebra,
is an associative $K$-algebra $R$ with two $K$-linear operators
$P,D:R\rightarrow R$ such that for any $x,y\in R$,

(I) (Rota-Baxter relation) $P(x)P(y)=P(xP(y))+P(P(x)y)+\lambda
P(xy);$

(II) ($\lambda$-differential relation) $ D(xy)=D(x)y+xD(y)+\lambda
D(x)D(y);$

(III) $D(P(x))=x$.

Hence, any  $\lambda$-differential Rota-Baxter
 algebra is an    associative algebra with  two linear operators $\Omega=\{P,
 D\}$.

Let  $K\langle X;\Omega \rangle$ be  the free  associative algebra
with multiple linear operators $\Omega$ on $X$, where
$\Omega=\{P,D\}$. For any $u\in \mathfrak{S}(X)$, $u$ has a unique
expression
$$
u=u_0P(b_1)u_1P(b_2)u_2\cdots P(b_n) u_n,
$$
where each $u_i\in (X\cup D(\mathfrak{S}(X)))^*$ and $b_i\in
\mathfrak{S}(X)$. Denote by
$$
wt(u)=(deg_{_{X}}(u),deg_{_{P}}(u),n, b_1, \cdots, b_n,
u_0,\cdots,u_n),
$$
where $deg_{_{P}}(u)$ is the number of $P$ in $u$. Also, for any
$u_t\in (X\cup D(\mathfrak{S}(X)))^* $, $u_t$ has a unique
expression
$$
u_t=u_{t_{1}}\cdots  u_{t_{k}}
$$
where each $u_{t_{j}}\in X\cup D(\mathfrak{S}(X))$.

 Let $X$  be well-ordered  and $u,v\in \mathfrak{S}(X)$. Order  $\mathfrak{S}(X)$ as follows:
\begin{equation}\label{o3}
u>v\Longleftrightarrow wt(u)>wt(v)\ \mbox{ lexicographically}
\end{equation}
 where for each $t$, $\ u_t>v_t$ if
 $$
 (deg_{_{X}}(u_t),deg_{_{P}}(u_t), u_{t_{1}},\cdots,u_{t_{k}})>
 (deg_{_{X}}(v_t),deg_{_{P}}(v_t), v_{t_{1}},\cdots,v_{t_{l}})\ \mbox{
 lexicographically}
 $$
where for each $j$,  $u_{t_{j}}>v_{t_{j}}$ if  one of the following
holds:

(a)  $u_{t_{j}},v_{t_{j}}\in X$ and $u_{t_{j}}>v_{t_{j}}$;

(b)  $u_{t_{j}}=D(u_{t_{j}}^{'}), v_{t_{j}}\in X$;

(c)  $u_{t_{j}}=D(u_{t_{j}}^{'}),v_{t_{j}}=D(v_{t_{j}}^{'})$ and
$u_{t_{j}}^{'}>v_{t_{j}}^{'}$.

\ \

Clearly, the order (\ref{o3}) is a monomial order on
$\mathfrak{S}(X)$.

Let $S$ be the set consisting of  the following
$\Omega$-polynomials:
\begin{enumerate}
\item[1.]
$P(x)P(y)-P(xP(y))-P(P(x)y)-\lambda P(xy)$,
\item[2.] $D(xy)-D(x)y-xD(y)-\lambda D(x)D(y)$,
\item[3.] $D(P(x))-x$,
\end{enumerate}
where $x,y\in \mathfrak{S}(X)$.

\begin{theorem}(\cite{BoChQiu-CD-Omega})\label{t2.2.4} With  the order
(\ref{o3}) on $\mathfrak{S}(X)$,  $S$ is a Gr\"{o}bner-Shirshov
basis in $K\langle X;\Omega\rangle$.
\end{theorem}

By Theorems \ref{CD-RotaB} and \ref{t2.2.4}, we obtain a normal form
$Irr(S)$ of  the free $\lambda$-differential Rota-Baxter algebra
$K\langle X;P,D|S\rangle$ which is a similar one in \cite{GK08}.

\subsection{Rota-Baxter algebras}
In the subsection, we consider Rota-Baxter algebras over a field of
characteristic 0.

The free non-commutative Rota-Baxter algebra is given by K.
Ebrahimi-Fard and L. Guo \cite{EG08a}. The free commutative
Rota-Baxter algebra is given by  G. Rota \cite{Ro69} and P. Cartier
\cite{Ca}.

Let $X$ be a nonempty set, $S(X)$ the free semigroup generated by
$X$ without identity and $P$ a symbol of a unary operation. For any
two nonempty sets $Y$ and $Z$, denote by
$$
{\Lambda}_{P}{(Y,Z)}=(\cup_{r\geq 0}(YP(Z))^{r}Y)\cup(\cup_{r\geq
1}(YP(Z))^{r})\cup(\cup_{r\geq 0}(P(Z)Y)^{r}P(Z))\cup(\cup_{r\geq
1}(P(Z)Y)^{r}),
$$
where for a set $T$,  $T^0$ means the empty set.

Define
\begin{eqnarray*}
\Phi_{0}&=&S(X)\\
\vdots\ \ & &\ \ \ \ \vdots\\
\Phi_{n}&=&{\Lambda}_{P}(\Phi_{0},\Phi_{n-1}) \\
\vdots\ \ & &\ \ \ \ \vdots\\
\end{eqnarray*}
Let
$$
\Phi(X)=\cup_{n\geq 0} \Phi_{n}.
$$
Clearly, $P(\Phi(X))\subset \Phi(X)$. If $u\in X\cup P({\Phi}(X))$,
then $u$ is called prime.  For any $u\in \Phi(X)$, $u$ has a unique
form
 $u=u_{1}u_{2}\cdots u_{n}$ where $u_{i}$  is prime, $i=1, 2,\dots,n$, and
 $u_i,u_{i+1}$ can not both have forms as $p(u'_i)$ and
 $p(u'_{i+1})$. If this is the case, then we define the breath of $u$ to be $n$, denoted by
 $bre(u)=n$.

For any $u\in \Phi(X)$ and for a set $T\subseteq X\cup \{P\}$,
denote by $deg_{T}(u)$ the number of occurrences of $t\in T$ in $u$.
 Let
$$
Deg(u)=(deg_{\{P\}\cup X}(u), deg_{\{P\}}(u)).
$$
We order $Deg(u)$ lexicographically.

Let $k\Phi(X)$ be a free  $k$-module with $k$-basis $\Phi(X)$ and
$\lambda \in k$ a fixed element. Extend linearly $P: \
k\Phi(X)\rightarrow k\Phi(X),\ u\mapsto P(u)$ where $u\in \Phi(X)$.

Now we define the multiplication in $k\Phi(X)$.

Firstly, for $u,v\in X\cup P(\Phi(X))$, define

$u\cdot v= \left\{ \begin{array}
              {l@{\quad}l}
              P(P(u')\cdot v')+P(u'\cdot P(v'))+\lambda P(u'\cdot v'), & \mbox{if}\ u=P(u'),v=P(v'); \\
              uv, & \mbox{otherwise}.
              \end{array} \right. $

Secondly, for any $u=u_{1}u_{2}\cdots u_{s}, v=v_{1}v_{2}\cdots
v_{l}\in \Phi(X)$ where $u_{i}, v_{j}$ are prime, $i=1,2,\dots,s,
j=1,2,\dots,l$, define
$
u\cdot v=u_{1}u_{2}\cdots u_{s-1}(u_{s}\cdot v_{1})v_{2}\cdots
v_{l}.
$

Equipping with the above concepts, $k\Phi(X)$ is the free
Rota-Baxter algebra with weight $\lambda$ generated by $X$ (see
\cite{EG08a}), denoted by $RB(X)$.

We have to order $\Phi(X)$. Let $X$ be a well-ordered set. Let us
define an order $>$ on $\Phi(X)$ by induction on the $Deg$-function.

For any $u,v\in \Phi(X)$, if $Deg(u)>Deg(v)$, then $u>v$.

If $Deg(u)=Deg(v)=(n,m)$, then we define $u>v$ by induction on
$(n,m)$.

If $(n,m)=(1,0),$ then $u,v\in X$ and we use the order on $X$.
Suppose that for $(n,m)$ the order is defined where
$(n,m)\geq(1,0)$. Let $(n,m)<(n',m')=Deg(u)= Deg(v)$.  If $u,v \in
P(\Phi(X))$, say $u=P(u')$ and $v=P(v')$, then $u>v$ if and only if
$u'>v'$ by induction. Otherwise $u=u_1u_2\cdots u_l$ and
$v=v_1v_2\cdots v_s$ where $l>1$ or $s>1$, then $u>v$ if and only if
$(u_1, u_2, \dots, u_l)>(v_1, v_2, \dots, v_s )$ lexicographically
by induction.

It is clear that $>$ is a monomial order on $\Phi(X)$. Throughout
this subsection, we will use this order.

It is noted that not each $S$-word is normal.

Let $f,g \in RB(X)$ be monic with $\overline{f}=u_{1}u_{2}\cdots
u_{n}$ where each $u_{i} $ is prime. Then, there are four kinds of
compositions.
\begin{description}
\item{(i)} If $u_{n}\in P(\Phi(X))$, then we define composition of
right multiplication  as\ \ $f\cdot u$ \ \ where $u\in P(\Phi(X)).$

\item{(ii)} If $u_{1}\in P(\Phi(X))$, then we define composition of
left multiplication as \ \ $u\cdot f$ \ \  where \ $u\in
P(\Phi(X)).$

\item{(iii)} If there exits a $w=\overline{f}a=b
\overline{g}$ \ where $fa$ is normal $f$-word and $bg$ is normal
$g$-word, $a, b \in \Phi(X)$ and $deg_{\{P\}\cup
X}(w)<deg_{\{P\}\cup X}(\overline{f})+deg_{\{P\}\cup
X}(\overline{g})$, then we define the intersection composition of
$f$ and $g$ with respect to $w$ as $(f,g)_{w}=f\cdot a-b\cdot g.$

\item{(iv)} If there exists a $w=\overline{f}=u|_{\overline{g}}$ where
$u\in \Phi^{\star}(X)$, then we define the inclusion composition of
$f$ and $g$ with respect to $w$ as $(f,g)_{w}=f-u|_{g}.$
\end{description}

\begin{theorem}{\bf(\cite{BoChDeng-RotaB}, CD-lemma for Rota-Baxter algebras)}\label{CD-RotaB}
Let $RB(X)$ be a free Rota-Baxter algebra over a field of
characteristic 0 and $S$  a set of monic polynomials in $RB(X)$ and
$>$ the monomial order on $\Phi(X)$ defined as before. Then the
following statements are equivalent.
\begin{enumerate}
\item[(i)] $S $ is a Gr\"{o}bner-Shirshov basis in $RB(X)$.
\item[(ii)] $  f\in Id(S)\Rightarrow
\bar{f}=u|_{\overline{s}}$  for some $u \in \Phi(X)$,\ $s\in S$.
\item[(iii)] $Irr(S) = \{ u\in \Phi(X) |  u \neq
v|_{\bar{s}}, s\in S, v|_{s}\ \mbox{is normal}\  s\mbox{-word}\}$ is
a $k$-basis of $RB(X|S)\\=RB(X)/Id(S)$.
\end{enumerate}
\end{theorem}

The following are some applications of the Theorem \ref{CD-RotaB}.

\begin{theorem}(\cite{BoChDeng-RotaB})\label{Deng9} Let $I$ be a well-ordered set,
$X=\{x_i|i\in I\}$ and the order on $\Phi(X)$ defined as before. Let
\begin{eqnarray}
\label{e4}&&f=x_{i}x_{j}-x_{j}x_{i}, \ \ i>j, \ i,j\in I\\
\label{e5}&&g=P(u)x_{i}-x_{i}P(u), \ \  u\in \Phi(X),\ i\in I
\end{eqnarray}
Let $S$ consist of (\ref{e4}) and (\ref{e5}). Then $S$ is a
Gr\"{o}bner-Shirshov basis in $RB(X)$.
\end{theorem}

By using Theorems \ref{CD-RotaB} and \ref{Deng9}, we get a normal
form $Irr(S)$ of the free  commutative Rota-Baxter algebra $RB(X|S)$
which is the same as one in \cite{GKg00b}.

\begin{theorem}(\cite{BoChDeng-RotaB})\label{Deng12}
Every countably generated Rota-Baxter algebra with weight 0 can be
embedded into a two-generated Rota-Baxter algebra.
\end{theorem}

An important application of Theorem \ref{CD-RotaB} is PBW theorem
for dendriform algebra which is a conjecture of L. Guo \cite{GL}.

A dendriform algebra (see \cite{Lo}) is a $k$-space $D$ with two
binary operations $\prec$ and $\succ$ such that for any $x,y,z\in
D$,
\begin{eqnarray*}
\nonumber&&(x\prec y)\prec z=x\prec(y\prec z+y\succ z)\\
&&(x\succ y)\prec z=x\succ(y\prec z)\\
\nonumber&&(x\prec y+x\succ y)\succ z=x\succ(y\succ z)
\end{eqnarray*}

Suppose that $(D, \prec,\succ)$ is a dendriform algebra over $k$
with a linear basis $X=\{x_i|i\in I\}$. Let $x_i\prec x_j=\{x_i\prec
x_j\}, x_i\succ x_j=\{x_i\succ x_j\}$, where $\{x_i\prec x_j\}$ and
$\{x_i\succ x_j\}$ are linear combinations of $x\in X$. Then $D$ has
an expression by generators and defining relations
$$
D=D(X|x_i\prec x_j=\{x_i\prec x_j\}, x_i\succ x_j=\{x_i\succ x_j\},
x_i,x_j\in X).
$$

Denote by
$$
U(D)=RB(X|x_iP(x_j)=\{x_i\prec x_j\}, P(x_i)x_j=\{x_i\succ x_j\},
x_i,x_j\in X).
$$
Then $U(D)$ is the universal enveloping Rota-Baxter algebra of $D$,
see \cite{EG08a}.

The following is the PBW theorem for dendriform algebras which is
proved by Yuqun Chen, Qiuhui Mo \cite{ChMoAMS}.

\begin{theorem}(\cite{ChMoAMS})
Every dendriform algebra over a field of characteristic 0 can be
embedded into its universal enveloping Rota-Baxter algebra.
\end{theorem}

\subsection{Right-symmetric (or pre-Lie)  algebras}

A non-associative $A$ is called a right-symmetric (or pre-Lie)
algebra if $A$ satisfies  the following identity $(x,y,z)=(x,z,y)$
for the associator $(x,y,z)=(xy)z-x(yz)$. It is a Lie admissible
algebra in a sense that $A^{(-)}=(A,[xy]=xy-yx)$ is a Lie algebra.

Let $X=\{x_i|i\in I \}$ be a set and for any $(u)\in X^{\ast \ast
}$, $|(u)|$ the length of the word $(u)$.

Let $I$ be a well-ordered set. We order $X^{**}$ by the induction on
the lengths of the words $(u)$ and $(v)$:
\begin{enumerate}
\item[(i)] \  If $|((u)(v))|=2$, then $(u)=x_i > (v)=x_j$ if and
only if $i>j$.
\item[(ii)] \ If $|((u)(v))|>2$, then $(u)>(v)$ if and
only if one of the following cases holds:
\begin{enumerate}
\item[(a)] $|(u)|>|(v)|$.
\item[(b)] \ If $|(u)|=|(v)|, \ (u)=((u_1)(u_2))$ and
$(v)=((v_1)(v_2))$, then $(u_1)>(v_1)$ or ($(u_1)=(v_1)$  and
$(u_2)>(v_2)$).
\end{enumerate}
\end{enumerate}

It is clear that the order $<$ on $X^{**}$ is well-ordered which is
called deg-lex order on non-associative words. We use this order
throughout this subsection.

We now cite the definition of good words (see \cite{Se94}) by
induction on  length:
\begin{enumerate}
\item[1)]$x$ is a good word for any $x\in X$.

Suppose that we define good words of length $<n$.

\item[2)] non-associative word $((v)(w))$ is called a good word if
\begin{enumerate}

\item[(a)] both $(v)$ and $(w)$ are good words,

\item[(b)] if $(v)=((v_1)(v_2))$, then $(v_2)\leq(w)$.
\end{enumerate}
\end{enumerate}

We denote $(u)$ by $[u]$, if $(u)$ is a good word.

Let $W$ be the set of all good words in the alphabet $X$ and
$RS\langle X\rangle$ the free right-symmetric algebra over a field
$k$ generated by $X$. Then $W$ forms a linear basis of the free
right-symmetric algebra $RS\langle X\rangle$, see \cite{Se94}.
Daniyar Kozybaev, Leonid Makar-Limanov, Ualbai Umirbaev
\cite{KMLU08} has proved that the deg-lex order on $W$ is monomial.

Let $S\subset RS\langle X\rangle$ be a set of monic polynomials and
$s\in S$. An $S$-word $(u)_s$ is called a normal $S$-word if
$(u)_{\bar s }=(a\bar s b)$ is a good word.

Let $f,g\in S$, $[w]\in W$ and $a,b\in X^{*}$. Then there are two
kinds of compositions.
\begin{enumerate}
\item[(i)] If $\bar{f}=[a\bar{g}b]$, then
$
(f,g)_{\bar{f}}=f-[agb]
$
is called composition of inclusion.

\item[(ii)] If $(\bar{f}[w])$ is not good, then
$
f\cdot [w]
$
is called composition of right multiplication.
\end{enumerate}

\begin{theorem}\label{t2.4.1}{\bf(\cite{BoChLi-GSB-rightsym}, CD-lemma
for right-symmetry algebras)}\ Let $S\subset RS\langle X\rangle$ be
a nonempty set of monic polynomials and the order $<$ be defined as
before.  Then the following statements are equivalent.
\begin{enumerate}
\item [(i)] $S$ is a Gr\"{o}bner-Shirshov basis in $RS\langle X\rangle$.

\item [(ii)] $f\in Id(S)\Rightarrow \bar f =[a\bar s b]$ for some $s\in S\
and\ a,b\in X^*$, where $[as b]$ is  a normal $S$-word.

\item [(iii)] $Irr(S)=\{[u]\in W |[u]\ne [a\bar s b]\ a,b\in X^*,\ s\in S \mbox{ and }
[as b] \mbox{ is a normal } S\mbox{-word}\}$ is a linear basis of
the algebra $RS\langle X | S\rangle=RS\langle X\rangle/Id(S)$.
\end{enumerate}
\end{theorem}

As an application, we have a GSB for universal enveloping
right-symmetric algebra of a Lie algebra.

\begin{theorem}\label{t2.4.2}(\cite{BoChLi-GSB-rightsym})\
Let $({\cal{L}},[,])$ be a Lie algebra with a well-ordered basis
$\{e_i|\ i\in I\}$. Let $ [e_i,e_j]=\sum\limits_{m}\alpha_{ij}^me_m
$, where $\alpha_{ij}^m\in k$. We denote
$\sum\limits_{m}\alpha_{ij}^me_m$  by $\{e_ie_j\}$. Let
$$
U({\cal{L}})=RS\langle \{e_i\}_I| \ e_ie_j-e_je_i=\{e_ie_j\}, \ i,j
\in I\rangle
$$
be the universal enveloping right-symmetric algebra of ${\cal{L}}$.
Let
\begin{eqnarray*}
S&=&\{f_{ij}=e_ie_j-e_je_i-\{e_ie_j\},\ i,j \in I \ \mbox{ and } \
i>j \}.
\end{eqnarray*}
Then the set $S$ is a Gr\"{o}bner-Shirshov basis in $RS\langle
X\rangle$ where $X=\{e_i\}_I$.
\end{theorem}

By Theorems \ref{t2.4.1} and \ref{t2.4.2}, we immediately have the
following PBW theorem for Lie algebra and right-symmetric algebra.

\begin{corollary}
(D. Segal \cite{Se94}) \ A Lie algebra ${\cal{L}}$ can be embedded
into its universal enveloping right-symmetric algebra $U({\cal{L}})$
as a subalgebra of $U({\cal{L}})^{(-)}$.
\end{corollary}

\subsection{$L$-algebras}

An $L$-algebra (see \cite{Ph}) is a $k$-space $L$ equipped with two
binary $k$-linear operations $\prec,~\succ : L{\otimes}L \rightarrow
L$ verifying the so-called entanglement relation:
\begin{equation*}
(x \succ y)\prec z = x \succ ( y \prec z), \ \ \forall  x, y ,z \in
L.
\end{equation*}

Let $\Omega=\{\succ,\ \prec\}=\Omega_2$ and
$\mathfrak{S}(X)=\bigcup_{n\geq0}\mathfrak{S}_{n}$ (see \S2.2). For
any $\delta\in \Omega,\ x,y\in X$, we write $\delta(x,y)=x\delta y$.

An $\Omega$-word $u$ is an $L$-word if $u$ is one of the following:
\begin{enumerate}
\item[i)]\  $u=x$, \ where $x\in X$.
\item[ii)]\ $u=v\succ w$, where $v$ and $w$ are $L$-words.
\item[iii)]\ $u=v\prec w$ with $v\neq v_1\succ v_2$, where $v_1,\ v_2,\ v,\ w$ are $L$-words.
\end{enumerate}

Let $L(X)$ be the free $L$-algebra generated by $X$. Then the set
$N$ of all $L$-words forms a normal form of $L(X)$, see
\cite{ChHuang-L-alg}, Corollary 3.3.

Let $X$ be a well-ordered set. Denote  $\succ$ by $\delta_1$ and
$\prec$ by $\delta_2$. Let $\delta_1<\delta_2$. For any $u \in
\mathfrak{S}(X)$, if $u=x\in X$, denote by
$
wt(u)=(1, x);
$
if $u=\delta_i(u_1,u_2)$ for some $u_1,~u_2\in \mathfrak{S}(X)$,
denote by
$
wt(u)=(deg_{_X}(u), \delta_i, u_1, u_2).
$
For any $u,v\in \mathfrak{S}(X)$, define
\begin{equation*}\label{1}
u>v\Longleftrightarrow wt(u)>wt(v) ~~~~~~~lexicographically
\end{equation*} by
induction on $|u|+|v|$.

It is clear that $>$ is a monomial order on $N$.

It is noted that not each $S$-word is normal.

In the paper \cite{ChHuang-L-alg},  CD-lemma for $\Omega$-algebras
is established, where $\Omega$ consists of $n$-ary operations,
$n\geq1$. This generalizes the result in V. Drensky and R. Holtkamp
\cite{Dren}. As a result,   the linear basis $N$ (the set of all
$L$-words) for the free $L$-algebra is obtained by using CD-lemma
for $\Omega$-algebras. Then  the CD-lemma for $L$-algebras  is
given. There are two kinds compositions: inclusion composition and
right multiplication composition. As applications, the following
embedding theorems for $L$-algebras are obtained:

\begin{theorem} (\cite{ChHuang-L-alg})\
1) Every countably generated $L$-algebra over a field $k$ can be
embedded into a two-generated $L$-algebra.

2) Every $L$-algebra over a field $k$ can be embedded into a simple
$L$-algebra.

3) Every countably generated $L$-algebra over a countable field $k$
can be embedded into a simple two-generated $L$-algebra.

4) Three arbitrary $L$-algebras $A$, $B$, $C$ over a field $k$ can
be embedded into a simple $L$-algebra generated by $B$ and $C$ if
$|k|\leq \dim(B*C)$ and $|A|\leq|B*C|$, where $B*C$ is the free
product of $B$ and $C$.
\end{theorem}

GSB of a free dialgebra and the free product of two $L$-algebras,
respectively are also given, and then the normal forms of such
algebras are obtained in \cite{ChHuang-L-alg}.

\subsection{Differential algebras}

Let ${\cal{D}}$ be a set of symbols,  ${A}$ an associative algebra
over $K$. Then $A$ is a differential algebra with differential
operators $\cal{D}$ or $\cal{D}$-algebra for short if for any
$\delta\in {\cal{D}}, \ a,b \in A$,
$$
\delta(ab)=\delta(a)\cdot b+a\cdot \delta(b).
$$

Let ${\cal{D}}=\{D_j|j\in J\}$ and $\mathbb{N}$ the set of
non-negative integers. For any
 $m\in \mathbb{N}$ and $\bar{j}=(j_1,\cdots,j_m)\in J^m$,
denote by $D^{\bar{j}}=D_{j_1}D_{j_2}\cdots D_{j_{m}}$ and
$D^{\omega}(X)=\{D^{\bar{j}}(x)|x\in X, \ \bar{j}\in J^m, \
m\in\mathbb{N}\}$, where $D^{\bar{j}}(x)=x$ if $\bar{j}\in J^0$. Let
$T=(D^{\omega}(X))^*$ be the free monoid generated by
$D^{\omega}(X)$.

Let ${\cal{D}}(X)$ be the free differential algebra with
differential operators $\cal{D}$ generated by $X$. Then $T$ is a
linear basis of ${\cal{D}}(X)$, see \cite{ChChLi-GSB-diff}.

Let  $X$, $J$ be well-ordered sets,
 $D^{\bar{i}}(x)=D_{i_1}D_{i_2}\cdots D_{i_{m}}(x)\in D^{\omega}(X)$  and
$$
wt(D^{\bar{i}}(x))=(x; m, i_1, i_2,\cdots, i_{m}).
$$
Then, we order $D^{\omega}(X)$ as follows:
$$
D^{\bar{i}}(x)> D^{\bar{j}}(y)\Longleftrightarrow
wt(D^{\bar{i}}(x))> wt(D^{\bar{j}}(y)) \ \mbox{ lexicographically}.
$$
It is easy to check this order is a well ordering on
$D^{\omega}(X)$.

Then deg-lex order on $T=(D^{\omega}(X))^*$ is monomial.

It is noted that each $S$-word is normal.

Let $f,g\in {\cal{D}}(X)$ be monic polynomials and $w,a,b\in T$.
Then there are two kinds of compositions.
\begin{enumerate}
\item[(i)] There are two sorts of composition of inclusion:

If $w=\bar{f}=a\cdot d^{\bar{j}}(\bar{g})\cdot b$, then the
composition is
$
(f,g)_{w}=f-a\cdot D^{\bar{j}}(g)\cdot b.
$
If $w=d^{\bar{i}}(\bar{f})=\bar{g}\cdot b$, then the composition is
$
(f,g)_{w}=D^{\bar{i}}(f)-g\cdot b.
$

\item[(ii)] Composition of intersection:

If $w=\bar{f}\cdot b=a\cdot d^{\bar{j}}(\bar{g})$ such that
$|\bar{f}|+|\bar{g}|>|w|$, then the composition is
$
(f,g)_{w}=f\cdot b-a\cdot D^{\bar{j}}(g).
$
In this case, we assume that $a,b\neq 1$.
\end{enumerate}

\begin{theorem}\label{t2.6.1}
\textbf{(\cite{ChChLi-GSB-diff}, CD-lemma for differential
algebras)} Let ${\cal{D}}(X)$ be the free differential algebra with
differential operators ${\cal{D}}=\{D_j|
 j\in J\}$,  $S\subset {\cal{D}}(X)$ a
monic subset and  $<$ the order on $T=(D^{\omega}(X))^*$ as before.
Then the following statements are equivalent.
\begin{enumerate}
\item[(i)]$S$ is a  Gr\"{o}bner-Shirshov basis in ${\cal{D}}(X)$.

\item[(ii)] $f\in Id(S)\Rightarrow\bar{f} =a\cdot
d^{\bar{i}}(\bar{s})\cdot b $ for some $s\in S,\ \bar{i}\in J^m, \
m\in \mathbb{N}$ and $a,b\in T.$

\item[(iii)] $Irr(S)=\{ u\in T \ | \ u\neq a\cdot
d^{\bar{i}}(\bar{s})\cdot b  \ for \ all \ s\in S, \  a,b\in T, \
\bar{i}\in J^m, \ m\in \mathbb{N}\}$ is a $k$-linear basis of
${\cal{D}}(X|S)={\cal{D}}(X)/Id(S)$.
\end{enumerate}
\end{theorem}

In the paper \cite{ChChLi-GSB-diff},  as applications of Theorem
\ref{t2.6.1}, there are given GSB for free Lie-differential algebras
and free commutative-differential algebras, respectively.

\subsection{Non-associative algebras  over a commutative algebra}

Let $K$ be a commutative associative $k-$algebra with unit, $X$ a
set and ${K}(X)$ the free non-associative algebra over $K$ generated
by $X$.

Let $[Y]$ denote the free abelian monoid generated by $Y$. Denote by
$$
N=[Y]X^{**}=\{u=u^{Y}u^{X}|u^{Y}\in [Y],u^{X}\in X^{**}\}.
$$
Then $N$ is a $k$-linear basis of the ``double" free non-associative
algebra ${k[Y]}(X)={k[Y]}\otimes k(X)$ over $k[Y]$ generated by $X$.

Suppose that both $>_{X}$ and $>_{Y}$ are monomial orders on
$X^{**}$ and $[Y]$, respectively. For any $u=u^Yu^X,v=v^Yv^X\in N$,
define
$$
u>v\Leftrightarrow %|u|>|v|  \ or \ (|u|=|v| \ and \ (
u^X>_{X}v^X \  or \ (u^X=v^X \ and \ u^Y>_{Y}v^Y).
$$
It is obvious that $>$ is a monomial order on $N$.

It is noted that each $S$-word is normal.

Let $f$ and $g$ be monic polynomials of $k[Y](X)$,  $w=w^Yw^X\in
[Y]X^{**}$ and $a,b,c\in X^*$, where $w^Y=L(\bar{f}^Y, \bar{g}^Y)=
L$ and $L(\bar{f}^Y, \bar{g}^Y)$ is the least common multiple of
$\bar{f}^Y$ and $\bar{g}^Y$ in  $k[Y]$. Then we have the following
two kinds compositions.

$(i)$\ $X$-inclusion

If $w^X=\bar{f}^X=(a(\bar{g}^X)b)$, then
$
(f,g)_{w}=\frac{L}{ \bar{f}^Y} f-\frac{L}{ \bar{g}^Y}(a(g)b)
$
is called the composition of $X$-inclusion.

$(ii)$\  $Y$-intersection only

If $|\bar{f}^Y|+|\bar{g}^Y|> |w^Y|$ and
$w^X=(a(\bar{f}^X)b(\bar{g}^X)c)$, then
$
(f,g)_{w}=\frac{L}{ \bar{f}^Y} (a(f)b(\bar{g}^X)c)-\frac{L}{
\bar{g}^Y} (a(\bar{f}^X)b(g)c)
$
is called the composition of $Y$-intersection only, where for $u\in
[Y],\ |u|$ means the degree of $u$.

A CD-lemma for non-associative algebras  over commutative algebras
is given (\cite{ChLiZeng-Cd-nonasso}, Theorem 2.6).  As an
application, the following embedding theorem is obtained.

\begin{theorem}(\cite{ChLiZeng-Cd-nonasso})\ Each countably generated
non-associative algebra over an arbitrary commutative algebra $K$
can be embedded into a two-generated non-associative algebra over
$K$.
\end{theorem}

\subsection{$n$-conformal algebras}

\begin{definition}
 Let $k$ be a field with characteristic $0$ and $C$ a vector
 space over $k$. Let $Z_+$ be the set of non-negative integers, $Z$ the integer ring and $n$
 a positive integer number.
 We associate to each $
 \overrightarrow{m}=(m_1,\cdots,m_n) \in Z_+^n$
 a bilinear product $\langle
 \overrightarrow{m}\rangle$ on $C$. Let
 $D_i:C\rightarrow C$ be linear mappings such that  $D_iD_j=D_jD_i, \ 1\leqslant i,j \leqslant n$.
 Then $C$ is an $n$-conformal algebra with derivations $D=\{D_1,\dots, D_n\}$ if the following axioms are
 satisfied:
\begin{enumerate}
\item[(i)]\ If $a,b \in C$, then there is an $\overrightarrow{N}(a,b)\in
Z_+^n$ such that
 $a \langle \overrightarrow{m}\rangle b=0$ if $\overrightarrow{m}\not\prec
 \overrightarrow{N}(a,b)$, where for any $\overrightarrow{m}=(m_1,\cdots,m_n),
\overrightarrow{l}=(l_1,\cdots,l_n)\in Z_+^n$,
\begin{eqnarray*}
&&\overrightarrow{m}=(m_1,\cdots,m_n)\prec
\overrightarrow{l}=(l_1,\cdots,l_n)\Leftrightarrow m_i\leq l_i, \
i=1,\dots,n \\
&& \mbox{and there exists }\ \ i_0, \ 1\leq i_0\leq n \ \ \mbox{
such that }\ \ \ m_{i_0}< l_{i_0}.
\end{eqnarray*}
 $\overrightarrow{N}:\ C\times C\rightarrow Z_+^n$ is called the locality
 function.
\item[(ii)]\ For any $a,b \in C, \ \overrightarrow{m} \in Z_+^n$,
      $D_i(a \langle \overrightarrow{m}\rangle b)=D_ia \langle \overrightarrow{m}\rangle b+
      a \langle \overrightarrow{m}\rangle D_ib,\ i=1,\dots,n.$
\item[(iii)]\ For any $a,b \in C, \ \overrightarrow{m} \in Z_+^n$,
    $D_ia \langle \overrightarrow{m}\rangle b=-m_ia \langle \overrightarrow{m}-\overrightarrow{e_i}\rangle
    b,\  i=1,\dots,n$, where
    $\overrightarrow{e_i}=(\underbrace{0,\cdots,0}_{i-1},1,0,\cdots,0)$.
\end{enumerate}
\end{definition}

For $\overrightarrow{m},\overrightarrow{s}\in Z_+^n$, put
$(-1)^{\overrightarrow{s}}=(-1)^{s_1+\cdots+s_n}$ and
$\binom{\overrightarrow{m}}{\overrightarrow{s}}=\binom{m_1}{s_1}\cdots\binom{m_n}{s_n}$.

\ \

An $n$-conformal algebra $C$ is associative if in addition the
associativity condition holds.

{\bf Associativity Condition:} For  any $a,b,c \in C$ and
$\overrightarrow{m},\overrightarrow{m}' \in Z_+^n$,
$$
   (a\langle \overrightarrow{m}\rangle b)\langle \overrightarrow{m}'\rangle
   c=\sum_{\overrightarrow{s}\in
   Z_+^n}(-1)^{\overrightarrow{s}}\binom{\overrightarrow{m}}{\overrightarrow{s}}a\langle
   \overrightarrow{m}-\overrightarrow{s}\rangle (b\langle
   \overrightarrow{m}'+\overrightarrow{s} \rangle c),
   $$
   or equivalently a $right$ analogy
   $$
   a\langle \overrightarrow{m}\rangle (b\langle \overrightarrow{m}'\rangle
   c)=\sum_{\overrightarrow{s}\in
   Z_+^n}(-1)^{\overrightarrow{s}}\binom{\overrightarrow{m}}{\overrightarrow{s}}(a\langle
   \overrightarrow{m}-\overrightarrow{s}\rangle b)\langle
   \overrightarrow{m}'+\overrightarrow{s}\rangle c.
   $$

\ \

We fix a locality function $\overrightarrow{N}:\ B\times
B\rightarrow Z_+^n$.

In the paper \cite{BoChZhang-n-conf}, by  a close analogy with the
paper by L. A. Bokut, Y. Fong, W.-F. Ke \cite{BFK04}, the free
associative $n$-conformal algebra
 is constructed and  CD-lemma for associative
$n$-conformal algebras is established. There are five kinds
compositions.

Let $C(B,\overrightarrow{N},D_1,\cdots,D_n)$ be the free associative
$n$-conformal algebra generated by $B$ with a fixed locality
function $\overrightarrow{N}$. Then the set of elements of the form
\begin{eqnarray*}\label{e3.1}
[u]=a_1\langle \overrightarrow{m}^{(1)}\rangle (a_2\langle
    \overrightarrow{m}^{(2)}\rangle(\cdots (a_k\langle
    \overrightarrow{m}^{(k)}\rangle D^{\overrightarrow{i}}a_{k+1})\cdots))
\end{eqnarray*}
forms a linear basis of $C(B,\overrightarrow{N},D_1,\cdots,D_n)$,
where $a_l\in B$, $\overrightarrow{m}^{(j)}\in Z_+^n$,
$\overrightarrow{m}^{(j)}\prec \overrightarrow{N}$, $1\leqslant l
\leqslant k+1, 1\leqslant j \leqslant k$,
ind(u)=$\overrightarrow{i}\in Z_+^n$, $|u|=k+1,k\geqslant 0$. We
shall refer to $[u]$ as D-free if ind(u)=$(0,\cdots,0)$, and we
shall say that $f$ is D-free if every normal word in $f$ is D-free.

\begin{theorem} {\bf(\cite{BoChZhang-n-conf}, CD-lemma for associative $n$-conformal
algebras)}\label{t2.8.2}\\
Let $S$ be a Gr\"{o}bner-Shirshov basis in
$C(B,\overrightarrow{N},D_1,\cdots,D_n)$. Then the set $Irr(S)$
forms a linear basis of the $n$-conformal algebra
$C(B,\overrightarrow{N},D_1,\cdots,D_n|S)$ with defining relations
$S$.
\par If $S$ is D-free, then the converse is true as well.
\end{theorem}

Remark that the condition ``the set $S$ is a Gr\"{o}bner-Shirshov
basis" is not equivalent to ``the set $Irr(S)$ forms a linear basis
of the $n$-conformal algebra
$C(B,\overrightarrow{N},D_1,\cdots,D_n|S)$".

As an application of Theorem \ref{t2.8.2},   there are constructed
GSB for loop Lie $n$-conformal algebras presented by generators and
defining relations in \cite{BoChZhang-n-conf}.

\subsection{$\lambda$-differential  associative algebras with
multiple operators}

A $\lambda$-differential associative algebra with multiple operators
is a $\lambda$-differential algebra $R$ with a set
 $\Omega$ of multi-linear  operators.

Free $\lambda$-differential associative algebra on $X$ with multiple
operators $\Omega$ is constructed in \cite{ChQiu-Cd-diff}.

There are two kinds of compositions: intersection and inclusion.

A CD-lemma for $\lambda$-differential  associative algebras with
multiple operators is established in \cite{ChQiu-Cd-diff}. As an
application, a Gr\"{o}bner-Shirshov basis for the free
$\lambda$-differential Rota-Baxter  algebra is given and then normal
forms is obtained for such an algebra which is the same as one in
\cite{BoChQiu-CD-Omega}.

\subsection{Categories}

Let ${\cal{C}}$ be a category. Let
\begin{eqnarray*}
k{\cal{C}}&=&\{f=\sum_{i=1}^n \alpha_i\mu_i|\alpha_i\in k, \
\mu_i\in mor({\cal{C}}),\ n\geq0, \\
&&\ \  \mu_{i}\ (0\leq i\leq n) \mbox{ have the same domains and the
same codomains} \}.
\end{eqnarray*}
Note that in $k{\cal{C}}$, for $f,\ g\in k{\cal{C}},\ f+g$ is
defined only if $f,\ g$ have the same domain and the same codomain.

A multiplication $\cdot $ in $k{\cal{C}}$ is defined by linearly
extending the usual compositions of morphisms of the category
${\cal{C}}$. Then ($k{\cal{C}},\ \cdot$) is called the category
partial algebra over $k$ relative to ${\cal{C}}$ and $kC(X)$ the
free category partial algebra generated by the graph $X$.

Just like the case of associative algebras, a CD-lemma for
categories is established, see \cite{BoChLi-Nankai}, Theorem 4.4. As
applications, they give GSB for the simplicial category and the
cyclic category respectively, cf. \cite{BoChLi-Nankai}, Theorems 5.1
and 5.4, see S. I. Gelfand, Y. I. Manin, \cite{Gelfand}, Chapter
3.1.

\subsection{Metabelian Lie algebras}

 Suppose that $L$ is a Lie algebra over a field.
As usual we set $L^{(0)}=L$, $L^{(n+1)}=[L^{(n)}, L^{(n)}]$.
 Then $L/L^{(2)}$ is called a metabelian Lie algebra. More precisely,
 the variety of metabelian Lie algebras is given by the identity
 $$
 (x_1x_2)(x_3x_4)=0.
 $$

Let us begin with the construction of a free metabelian Lie algebra.
Let $X$ be a linear ordered set, $Lie(X)$  the free Lie algebra
generated by $X$. Then $L_{(2)}(X)=Lie(X)/Lie(X)^{(2)}$ is the free
metabelian Lie algebra generated by $X=\{a_i\}$. We call a monomial
is left-normed if it is of the form $(\cdots((ab)c)\cdots)d$. In the
sequel, the left-normed brackets in the expression is omitted
without confusion. For an arbitrary set of indices
$j_1,j_2,\cdots,j_m$, define
$$
\langle a_{j_1}\cdots a_{j_m}\rangle=a_{i_1}\cdots a_{i_m},
$$
where $ a_{i_1}\leq \cdots \leq a_{i_m} $ and $i_1,i_2,\cdots,i_m$
is a permutation of the indices $j_1,j_2,\cdots,j_m$.

Let
$$
R=\{u=a_0a_1a_2\cdots a_n \ | \ u  \ \mbox{is left-normed}, a_i\in
X,  \ a_0>a_1\leq a_2\leq \cdots \leq a_n ,
 n\geq1\}
$$
and $N=X\cup R$. Then $N$ forms a linear basis of the free
metabelian Lie algebra $L_{(2)}(X)$. Therefore, for any $f\in
L_{(2)}(X)$, $f$ has a presentation $f=f^{(1)}+f^{(0)}$, where
$f^{(1)}\in kR$ and $f^{(0)}\in kX$.

Moreover, the multiplication of the elements of $N$ is the
following, $u\cdot v=0$ if both $u,v\in R$, and
$$
a_0a_1a_2\cdots a_n\cdot b=\left\{
\begin{array}{l}
a_0\langle a_1a_2\cdots a_nb\rangle \ \ \ \ \ \ \ \ \ \ \ \ \ \ \ \ \ \ \ \ \ \  \ \ \  \mbox{if} \ a_1\leq b ,\\
a_0ba_1a_2\cdots a_n-a_1b\langle a_0a_2\cdots a_n\rangle \ \ \ \
\mbox{if}  \ a_1> b .
\end{array}\right.
$$

If $u=a_0a_1\cdots a_n\in R$, then the words $a_i$, $a_0\langle
a_{i_1}\cdots a_{i_s}\rangle$ ($s\leq n, \ a_{i_1}\cdots a_{i_s}$ is
a subsequence of the sequence $a_1\cdots a_n$) are called subwords
of the word $u$.

Define the length of monomials in $N$:
$$
|a|=1,  \ |a_0a_1a_2\cdots a_n|=n+1.
$$
Now we order the set $N$ degree-lexicographically.

Then we have  $\overline{a_0a_1a_2\cdots a_n\cdot b}=a_0\langle
a_1a_2\cdots a_nb\rangle$ and $|u\cdot b|=|u|+1$. For any $f\in
L_{(2)}(X)$, we call $f$ to be monic, $(1)$-monic and $(0)$-monic if
the coefficients of $\bar{f}$, $\overline{f^{(1)}}$ and
$\overline{f^{(0)}}$ is $1$ respectively.

Let $S\subset L_{(2)}(X)$. Then the following two kinds of
polynomial are called normal $S$-words:
\begin{enumerate}
\item[(i)] $sa_1a_2\cdots a_n$, where  $a_1\leq a_2\leq
\cdots \leq a_n$, $s\in S$, $\bar{s}\neq a_1$ and $n\geq0$;
\item[(ii)] $us$, where $u\in R$, $s\in S$ and $\bar{s}\neq u$.
\end{enumerate}

Let $f$ and $g$ be momic polynomials of $L_{(2)}(X)$. We define
seven different types of compositions as follows:
\begin{enumerate}
\item[1.] If $\bar{f}=a_0a_1\cdots a_n$, $\bar{g}=a_0b_1\cdots b_m$,
$(n,m\geq0)$, then let $w=a_0\langle lcm(ab) \rangle$, where $lcm(ab)$
denotes the least common multiple of associative words $a_1\cdots a_n$
and $b_1\cdots b_m$. The composition of type $I$
of $f$ and $g$ relative to $w$ is defined by
$$
C_I( f,g)_{w}=f\langle\frac{lcm(ab)}{a_1\cdots
a_n}\rangle-g\langle\frac{lcm(ab)}{b_1\cdots b_m}\rangle.
$$

\item[2.] If $\bar{f}=\overline{f^{(1)}}=a_0a_1\cdots a_n$,
$\overline{g^{(0)}}=a_i$ for some $i\geq2$ or
$\overline{g^{(0)}}=a_1$ and $a_0>a_2$, then let $w=\bar{f}$ and the
composition of type $II$
of $f$ and $g$ relative to $w$ is defined by
$$
C_{II}( f,g)_{w}=f-\alpha^{-1}a_0a_1\cdots \hat{a_i}\cdots a_n\cdot
g,
$$
where $\alpha$ is the coefficient of $\overline{g^{(0)}}$.

\item[3.] If $\bar{f}=\overline{f^{(1)}}=a_0a_1\cdots a_n$, $\bar{g}=
\overline{g^{(0)}}=a_1$ and $a_0\leq a_2$ or $n=1$, then let $w=\bar{f}$
and the  composition of type $III$
of $f$ and $g$ relative to $w$ is defined by
$$
C_{III}( f,g)_{\bar{f}}=f+ga_0a_2\cdots a_n.
$$

\item[4.] If $\bar{f}=\overline{f^{(1)}}=a_0a_1\cdots a_n$, $g^{(1)}\neq0$,
$\overline{g^{(0)}}=a_1$ and $a_0\leq a_2$ or $n=1$, then for any $a<a_0$
and $w=a_0\langle a_1\cdots a_na\rangle$, the  composition of type $IV$
of $f$ and $g$ relative to $w$ is defined by
$$
C_{IV}( f,g)_{w}=fa-\alpha^{-1}a_0aa_2\cdots a_n\cdot g,
$$
where $\alpha$ is the coefficient of $\overline{g^{(0)}}$.

\item[5.] If $\bar{f}=\overline{f^{(1)}}=a_0a_1\cdots a_n$, $g^{(1)}\neq0$
and $\overline{g^{(0)}}=b\notin \{a_i\}_{i=1}^n$, then let $w=
a_o\langle a_1\cdots a_nb\rangle$ and the  composition of type $V$
of $f$ and $g$ relative to $w$ is defined by
$$
C_{V}( f,g)_{w}=fb-a_0a_1\cdots a_n\cdot g,
$$
where $\alpha$ is the coefficient of $\overline{g^{(0)}}$.

\item[6.] If $\overline{f^{(0)}}=\overline{g^{(0)}}=a$ and $f^{(1)}\neq0$,
then for any $a_0a_1\in R$ and $w=a_0\langle a_1a\rangle$, the  composition of type $VI$
of $f$ and $g$ relative to $w$ is defined by
$$
C_{VI}( f,g)_{w}=(a_0a_1)(\alpha^{-1}f-\beta^{-1}g),
$$
where $\alpha$ and $\beta$ are the coefficients of
$\overline{f^{(0)}}$ and $\overline{g^{(0)}}$ respectively.

\item[7.] If $f^{(1)}\neq0$, $g^{(1)}\neq0$ and $\overline{f^{(0)}}=
a>\overline{g^{(0)}}=b$, then for any $a_0>a$ and $w=a_0ba$, the  composition of type $VII$
of $f$ and $g$ relative to $w$ is defined by
$$
C_{VII}( f,g)_{w}=\alpha^{-1}(a_0b)f-\beta^{-1}(a_0a)g,
$$
where $\alpha$ and $\beta$ are the coefficients of
$\overline{f^{(0)}}$ and $\overline{g^{(0)}}$ respectively.
\end{enumerate}

Then, in \cite{ChCh-Matabel}, CD-lemma for metabelian Lie algebras
is given.

\ \

Remark: In V. V. Talapov \cite{Talapov}, there is  a mistake for
compositions that is now modified.

\subsection{$S$-act algebras}

Let  $A$ be an associative algebra over $k$ and $S$ a monoid of
linear operators on $A$. Then $A$ is called an $S$-act algebra if
$A$ is an $S$-act with the action $s(a)$ satisfying
$$
s(ab)=s(a)s(b),\ \ s\in S, \ a,b\in A.
$$

In the paper \cite{Xia},  the ``double free" $S$-act algebra (i.e.,
a free $S$-act algebra, where $S$ is a free semigroup) is
constructed. Then  Gr\"{o}bner-Shirshov bases theory for $S$-act
algebras is established, where $S$ is an arbitrary semigroup.  As an
application, a Gr\"{o}bner-Shirshov basis of free Chinese monoid-act
algebra is given and hence  a linear basis of free Chinese
monoid-act algebra is obtained.

\subsection{Operads}

We present elements of the free operad by trees.

Let $\mathscr{V}=\bigcup_{n=1}^{\infty}\mathscr{V}_{n},$ where
$\mathscr{V}_{n}=\{\delta_i^{(n)}|i\in I_n\}$ is the set of $n$-ary
operations.

A tree with $n$ leaves is called decorated if we  label the leaves
by $[n]=\{1,2,3,\ldots,n\}, n\in \mathbb{N}$ and each vertex by an
element in $\mathscr{V}$.

A decorated tree is called tree monomial if for each vertex, the
minimal value on the leaves of the left subtree is always less than
that of the right subtree.

For example, $\delta_1^{(2)}(\delta_2^{(2)}(1,3),2),\
\delta_1^{(2)}(1,\delta_2^{(2)}(2,3))$  are  tree monomials, but
$\delta_1^{(2)}(1,\delta_2^{(2)}(3,2))$ is not a tree monomial.

Let $\mathscr{F}_\mathscr{V}(n)$ be the set of all tree monomials
with $n$ leaves and $T=\cup_{n\geq1}\mathscr{F}_\mathscr{V}(n)$. For
any $\alpha=\alpha(x_1,\dots, x_n)\in \mathscr{F}_\mathscr{V}(n),\
\beta\in \mathscr{F}_\mathscr{V}(m)$, define the shuffle composition
$\alpha \circ_{i,\sigma} \beta$ as follows
$$
\alpha(x_1,
\ldots,x_{i-1},\beta(x_i,x_{\sigma(i+1)},\ldots,x_{\sigma(i+m-1)}),x_{\sigma(i+m)},\ldots,x_{\sigma(m+n-1)})
$$
which is in $\mathscr{F}_\mathscr{V}(n+m-1)$, where $1\leq i\leq n$
and the bijection
$\sigma:\{i+1,\ldots,m+n-1\}\rightarrow\{i+1,\ldots,m+n-1\}$ is an
$(m-1,n-i)$-shuffle; that is
\begin{eqnarray*}
&&\sigma(i+1)<\sigma(i+2)<\dots<\sigma(i+m-1),\\
&&\sigma(i+m)<\sigma(i+m+1)<\dots<\sigma(n+m-1).
\end{eqnarray*}
For example,

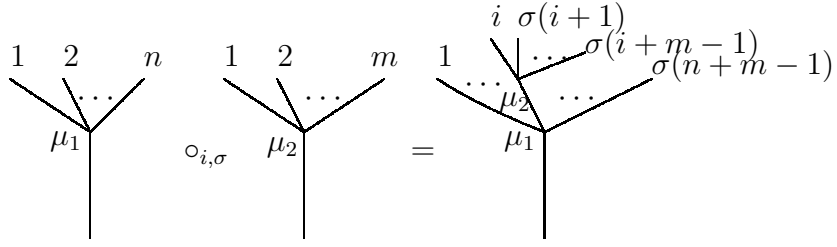
\begin{picture}(400,110)
\put(135,50){$\circ_{i, \sigma}$}\put(220,50){$=$}
\put(70,85){1}\put(90,85){2}\put(95,70){$\cdots$} \put(120,85){$n$}
\put(85,55){$\mu_1$}\put(150,85){1}\put(170,85){2}
\put(180,70){$\cdots$}\put(205,85){$m$}\put(166,53){$\mu_2$}
\put(230,85){1}\put(240,75){$\cdots$}\put(253,70){$\mu_2$}
\put(250,100){$i$}\put(260,100){$\sigma(i+1)$}
\put(265,85){$\cdots$}
\put(285,90){$\sigma(i+m-1)$}\put(275,70){$\cdots$}
\put(310,82){$\sigma(n+m-1)$}\put(255,55){$\mu_1$}

\qbezier(70,80)(85,70)(100,60)\qbezier(90,80)(95,70)(100,60)
\qbezier(120,80)(110,70)(100,60) \qbezier(100,60)(100,40)(100,20)

\qbezier(150,80)(165,70)(180,60)\qbezier(170,80)(175,70)(180,60)
\qbezier(210,80)(195,70)(180,60)\qbezier(180,60)(180,40)(180,20)

\qbezier(230,80)(245,70)(270,60)\qbezier(260,80)(265,70)(270,60)
\qbezier(250,95)(255,87.5)(260,80)\qbezier(260,95)(260,87.5)(260,80)
\qbezier(285,90)(272.5,85)(260,80)\qbezier(310,80)(290,70)(270,60)
\qbezier(270,60)(270,40)(270,20)
\end{picture}

Then with the shuffle composition, the set $T$ is freely generated
by $\mathscr{V}$.

Let $\mathscr{F}_\mathscr{V}=kT$ be  a $k$-space with $k$-basis $T$.
Then with the shuffle compositions $\circ_{i,\sigma}$,
$\mathscr{F}_\mathscr{V}$ is called the free shuffle operad.

Let $S$ be a homogeneous subset of $\mathscr{F}_\mathscr{V},\ s\in S$. The $S$-word $%
u|_s$ is defined as in \S1 (in this case, the $\star$ depends on
$s$). That means
\begin{enumerate}
\item[(i)] $s$ is an $S$-word.
\item[(ii)] If $u|_s$ is an $S$-word, then for any $v\in T$,
$u|_s\circ_{i,\sigma} v$ and $v\circ_{j,\tau} u|_s$ are both
$S$-words if the shuffle compositions are defined.
\end{enumerate}

It is clear that the ideal $Id(S)$ of $\mathscr{F}_\mathscr{V}$ is
the set of linear combination of $S$-words.

A well ordering $>$ on $T$ is called monomial (admissible) in the
sense that
$$
\alpha>\beta\Rightarrow u|_\alpha>u|_\beta \  \mbox{ for any } u\in
T.
$$

Suppose that $T$ equips with a monomial order. Then each $S$-word is
a normal $S$-word.

 For example, the following
order $>$ on $T$ is monomial, see \cite{DK10}, Proposition 5.

For any $\alpha=\alpha(x_1,\dots, x_n)\in
\mathscr{F}_\mathscr{V}(n),\ \alpha$ has an unique expression
$$
\alpha=(path(1),\dots, path(n),x_1,\dots,x_n),
$$
where each $path(i)\in\mathscr{V}^*$ is the unique path from the
root to the leaf $i$. If this is the case, we denote by
$$
wt(\alpha)=(n,path(1),\dots, path(n),x_1,\dots,x_n).
$$

Let $\mathscr{V}$ be a well-ordered set. We order $\mathscr{V}^*$ by
deg-lex order. We order numbers by inverse order, i.e., $x>y$ if and
only if as natural number, $x$ is less than $y$.

Now, for any $\alpha,\beta\in T$, we define
$$
\alpha>\beta \Leftrightarrow wt(\alpha)>wt(\beta)\ \ \
lexicographically.
$$

 For  $a,b \in T$, we say $a$ is divisible by $b$, if there
exists a subtree of the underlying tree of $a$ for which the
corresponding tree monomial $a'$ is equal to $b$.

An element of $\mathscr{F}_\mathscr{V}$ is said to be homogeneous if
all  tree monomials that occur in this element with nonzero
coefficients have the same arity degree (but not necessarily the
same operation degree).

\begin{definition}
Let $f, g$ be two monic homogeneous elements of
$\mathscr{F}_\mathscr{V}$. We define two types of compositions as
follows:
\begin{enumerate}
\item[1.] If $w=\bar{f}=u|_{\bar{g}}$ for some $g$-word $u|_{{g}} $, then
$ (f,g)_w=f-u|_{{g}} $ is called the inclusion composition  of $f$
and $g$.

\item[2.] If $w=a\circ_{i,\sigma}\bar{f}=\bar{g}\circ_{j,\tau} b$
(or $w=\bar{f}\circ_{i,\sigma}a=\bar{g}\circ_{j,\tau} b$) for some
$a,b\in T$ and the number of vertices of the underlying tree of $w$
is less than the total number of vertices for $f$ and $g$, then $
(f, g)_w=a\circ_{i,\sigma} f-g\circ_{j,\tau} b $ ($(f,
g)_w=f\circ_{i,\sigma} a-g\circ_{j,\tau} b $) is called the
intersection composition of $f$ and $g$.
\end{enumerate}
\end{definition}

\begin{theorem}{\bf(\cite{DK10}, Theorem 1, CD-lemma
for shuffle operads)}\ Let the notations be as above.
$S\subset\mathscr{F}_\mathscr{V}$  the nonempty set of monic
homogeneous elements and $<$ a monomial order on $T$. Then the
following statements are equivalent.
 \begin{enumerate}
\item[(i)] $S$ is a Gr\"{o}bner-Shirshov basis in
$\mathscr{F}_\mathscr{V}$.

\item[(ii)] $f\in Id(S) \Rightarrow \bar{f}=u|_{\bar{s}}$  for some $S$-word
$u|_s$.

\item[(iii)] $Irr(S)=\{u\in T| u \neq
v|_{\bar{s}}
 \mbox{ for  any } \ S\mbox{-word }\ v|_s\}$
is a $k$-basis of $\mathscr{F}_\mathscr{V}/Id(S)$.
\end{enumerate}
\end{theorem}

As applications, in \cite{DK10}, the authors compute the
Gr\"{o}bner-Shirshov bases for some well-known operads, such as the
operad Lie of Lie algebras, the operad As of associative algebras,
the operad PreLie of pre-Lie algebras and so on.

In the paper of V. Dotsenko, M. V. Johansson \cite{DJo}, an
implementation of the algorithm for computing Gr\"{o}bner-Shirshov
bases for operads is given.

\section{ Applications of known CD-lemmas}

\subsection{Artin-Markov normal form for braid
groups in Artin-Burau generators}

Let $B_{n}$ denote the braid group of type  ${\bf \bf{A}_{n}}$.
Then
$$
B_{n}=gp\langle \sigma_1,\dots, \sigma_n\ |\
\sigma_j\sigma_i=\sigma_i\sigma_j\ (j-1>i),\
\sigma_{i+1}\sigma_i\sigma_{i+1}=\sigma_i\sigma_{i+1}\sigma_i,\
1\leq i\leq n-1\rangle.
$$

Let $X=Y\dot{\cup}Z$,  $Y^*$ and  $Z$  be well-ordered. Suppose that
the order on $Y^*$ is monomial. Then, any word in $X$ has the form
$u=u_0z_1u_1\cdots z_{k}u_{k}$, where $k\geq 0,\ u_i\in Y^*,\
z_{i}\in Z$. Define the inverse weight of the word $u\in
 X^*$ by
$$
inwt(u)=( k, u_{k}, z_{k}, \cdots,u_1, z_{1}, u_{0} ).
$$
Now we order the inverse weights lexicographically as follows
$$
u>v\Leftrightarrow inwt(u)>inwt(v).
$$
Then we call the above order the inverse tower order. Clearly, this
order is  a monomial order on $X^*$.

In case $Y=T\dot{\cup} U$ and $Y^*$ is endowed with the inverse
tower order, we call the order of words in $X$ the inverse tower
order of words relative to the presentation
$$
X=(T\dot{\cup} U)\dot{\cup} Z.
$$
In general, we can define the inverse tower order of $X$-words
relative to the presentation
$$
X=(\cdots
(X^{(n)}\dot{\cup}X^{(n-1)})\dot{\cup}\cdots)\dot{\cup}X^{(0)},
$$
where $X^{(n)}$-words are endowed by a monomial order.

In the braid group $B_{n}$, we now introduce a new set of generators
which are called  the Artin-Burau generators. We set
\begin{eqnarray*}
&&s_{i,i+1}=\sigma_{i}^{2},\ \ \
s_{i,j+1}=\sigma_{j}\cdots\sigma_{i+1}\sigma_{i}^{2}\sigma_{i+1}^{-1}\cdots\sigma_{j}^{-1},
\ \ \ \  1\leq i<j\leq n-1;\\
&&\sigma_{i,j+1}=\sigma_{i}^{-1}\cdots\sigma_{j}^{-1},\ \ \ 1\leq
i\leq j\leq n-1;\ \ \ \sigma_{ii}=1,\ \ \ \{a,b\}=b^{-1}ab.
\end{eqnarray*}
Form the set
$$
S_{j}=\{s_{i,j},s_{i,j}^{-1}, \ 1\leq i,j<n\} \ \mbox{ and } \
\Sigma^{-1}=\{\sigma_{1}^{-1},\cdots\sigma_{n-1}^{-1}\}.
$$
Then the set
$$
S=S_{n}\cup S_{n-1}\cup \cdots \cup S_{2}\cup\Sigma^{-1}
$$
generates $B_{n}$ as a semigroup.

Now we order the set S in the following way:
$$
S_{n}<S_{n-1}<\cdots < S_{2}<\Sigma^{-1},
$$
and
$$
s_{1,j}^{-1}< s_{1,j}< s_{2,j}^{-1}<\cdots<s_{j-1,j} \ , \ \ \
\sigma_{1}^{-1}<\sigma_{2}^{-1}<\cdots\sigma_{n-1}^{-1}.
$$
With the above notation, we now order the $S$-words by using the
inverse tower order, according to the fixed presentation of $S$ as
the union of $S_{j}$ and $\Sigma^{-1}$. We order the $S_{n}$-words
by the $deg-inlex$ order, i.e., we first compare the words by length
and then by inverse lexicographical order, starting from their last
letters.

\begin{lemma}
The following Artin-Markov relations hold in the braid group
$B_{n}$. For $\delta=\pm1$,
\begin{eqnarray}\label{e77}
&&\sigma_{k}^{-1}s_{i,j}^{\delta}=s_{i,j}^{\delta}\sigma_{k}^{-1},
k\neq i-1,i,j-1,j\\
&&\sigma_{i}^{-1}s_{i,i+1}^{\delta}=s_{i,i+1}^{\delta}\sigma_{1}^{-1}\\
&&\sigma_{i-1}^{-1}s_{i,j}^{\delta}=s_{i-1,j}^{\delta}\sigma_{i-1}^{-1}\\
&&\sigma_{i}^{-1}s_{i,j}^{\delta}=\{s_{i+1,j}^{\delta},s_{i,i+1}\}\sigma_{i}^{-1}\\
&&\sigma_{j-1}^{-1}s_{i,j}^{\delta}=s_{i,j-1}^{\delta}\sigma_{j-1}^{-1}\\
&&\sigma_{j}^{-1}s_{i,j}^{\delta}=\{s_{i,j+1}^{\delta},s_{j,j+1}\}\sigma_{j}^{-1}
\end{eqnarray}
for $\ i<j<k<l,\ \varepsilon = \pm1 $,
\begin{eqnarray}
&&s_{j,k}^{-1}s_{k,l}^{\varepsilon}=\{s_{k,l}^{\varepsilon},s_{j,l}^{-1}\}s_{j,k}^{-1}\\
&&s_{j,k}s_{k,l}^{\varepsilon}=\{s_{k,l}^{\varepsilon},s_{j,l}s_{k,l}\}s_{j,k}\\
&&s_{j,k}^{-1}s_{j,l}^{\varepsilon}=\{s_{j,l}^{\varepsilon},s_{k,l}^{-1}s_{j,l}^{-1}\}s_{j,k}^{-1}\\
&&s_{j,k}s_{j,l}^{\varepsilon}=\{s_{j,l}^{\varepsilon},s_{k,l}\}s_{j,k}\\
&&s_{i,k}^{-1}s_{j,l}^{\varepsilon}=\{s_{j,l}^{\varepsilon},s_{k,l}s_{i,l}s_{k,l}^{-1}s_{i,l}^{-1}\}s_{i,k}^{-1}\\
&&s_{i,k}s_{j,l}^{\varepsilon}=\{s_{j,l}^{\varepsilon},s_{i,l}^{-1}s_{k,l}^{-1}s_{i,l}s_{k,l}\}s_{i,k}
\end{eqnarray}
for $j<i<k<l$ or $ i<k<j<l, \ and  \ \varepsilon, \ \delta=\pm 1 $,
\begin{eqnarray}
&&s_{i,k}^{\delta}s_{j,l}^{\varepsilon}=s_{j,l}^{\varepsilon}s_{i,k}^{\delta}
\end{eqnarray}
and
\begin{eqnarray}
&&\sigma_{j}^{-1}\sigma_{k}^{-1}=\sigma_{k}^{-1}\sigma_{j}^{-1}, \  \ j<k-1\\
&&\sigma_{j,j+1}\sigma_{k,j+1}=\sigma_{k,j+1}\sigma_{j-1,j} \ ,\ \ j<k\\
&&\sigma_{i}^{-2}=s_{i,i+l}^{-1}\\
\label{e87}&&s_{i,j}^{\pm1}s_{i,j}^{\mp1}=1
\end{eqnarray}
\end{lemma}

\begin{theorem}(\cite{BoChaSh07})
The Artin-Markov relations (\ref{e77})-(\ref{e87}) form a
Gr\"{o}bner-Shirshov basis of the braid group $B_{n}$ in terms of
the Artin-Burau generators relative to the inverse tower order of
words.
\end{theorem}

In the paper \cite{BoChaSh07}, it was claimed that some compositions
are trivial. The paper \cite{ChMo-GSB-Braid-Artin} supported the
claim and all compositions worked out explicitly.

\subsection{Braid groups in Artin-Garside generators}

The Artin-Garside generators of the braid group $B_{n+1}$ are
$\sigma_i,\ 1\leq i\leq n,\ \bigtriangleup,\ \bigtriangleup^{-1}$
(Garside 1969), where $\bigtriangleup=\Lambda_1\cdots\Lambda_n$,
with $\Lambda_i=\sigma_1\cdots\sigma_i$.

Let us order $\bigtriangleup^{-1}<\bigtriangleup<
\sigma_1<\dots<\sigma_n$. We order
$\{\bigtriangleup^{-1},\bigtriangleup,\sigma_1,\dots,\sigma_n\}^*$
by deg-lex order.

By $V(j, i),W(j, i), \dots$, where $j\leq i$, we understand positive
words in the letters $\sigma_j,\ \sigma_{j+1},\dots,\sigma_i$. Also
$V(i + 1,i) = 1,\ W(i + 1,i) = 1,\dots$.

Given $V = V(1, i )$, let $V^{(k)},\ 1\leq k\leq n-i$ be the result
of shifting in $V$ all indices of all letters by $k,\
\sigma_1\mapsto \sigma_{k+1},\dots, \sigma_i\mapsto\sigma_{k+i}$,
and we also use the notation $V^{(1)}=V'$. We write
$\sigma_{ij}=\sigma_i\sigma_{i-1}\dots\sigma_{j},\ j\leq i-1,\
\sigma_{ii}=\sigma_i,\ \sigma_{ii+1}=1$.

\begin{theorem}(\cite{Bo08})\
 A Gr\"{o}bner-Shirshov basis of $B_{n+1}$ in
the Artin-Garside generators consists of the following relations:
\begin{eqnarray*}
&&\sigma_{i+1}\sigma_{i}V(1,i-1)W(j, i )\sigma_{i+1j} =
\sigma_{i}\sigma_{i+1}\sigma_{i}V(1,i-1)\sigma_{ij}W( j, i )',\\
&&\sigma_{s}\sigma_{k}=\sigma_{k}\sigma_{s},\ s-k\geq 2,\\
&&\sigma_{1}V_1\sigma_{2}\sigma_{1}V_2\cdots
V_{n-1}\sigma_{n}\cdots\sigma_{1}= \bigtriangleup
V_1^{(n-1)}V_2^{(n-2)}\cdots V_{(n-1)}',\\
&&\sigma_l\bigtriangleup^{-1}=\bigtriangleup^{-1}\sigma_{n-l+1},\
1\leq l\leq n,\\
&&\bigtriangleup\bigtriangleup^{-1}=1,\
\bigtriangleup^{-1}\bigtriangleup=1,
\end{eqnarray*}
where $1\leq i \leq n-1,\ 1\leq j \leq i+1,\ W$ begins with
$\sigma_i$ if it is not empty, and $V_i = V_i (1, i )$.
\end{theorem}

As results, we have the following corollaries.

\begin{corollary}
The $S$-irreducible normal form of each word of $B_{n+1}$ coincides
with the Garside normal form of the word.
\end{corollary}

\begin{corollary}
(Garside (1969))\ The semigroup of positive braids $B^+_{n+1}$ can
be embedded into a group.
\end{corollary}

\subsection{Braid groups in Birman-Ko-Lee generators}

Recall that the Birman-Ko-Lee generators $\sigma_{ts}$ of the braid
group $B_n$ are the elements
$$
\sigma_{ts}=(\sigma_{t-1}\sigma_{t-2}\dots\sigma_{s+1})\sigma_{s}
(\sigma_{s+1}^{-1}\cdots\sigma_{t-2}^{-1})\sigma_{t-1}^{-1}.
$$
Then  $B_n$ has an expression.
\begin{eqnarray*}
B_n&=&gp\langle \sigma_{ts},\ n\geq
t>s\geq1|\sigma_{ts}\sigma_{rq}=\sigma_{rq}\sigma_{ts},\
(t-r)(t-q)(s-r)(s-q)>0,\\
&&\ \ \
\sigma_{ts}\sigma_{sr}=\sigma_{tr}\sigma_{ts}=\sigma_{sr}\sigma_{tr},
n\geq t>s>r\geq1\rangle.
\end{eqnarray*}
Denote by $\delta=\sigma_{nn-1}\sigma_{n-1n-2}\cdots\sigma_{21}$.

Let us order $\delta^{-1}<\delta<\sigma_{ts}<\sigma_{rq} \ \mbox{
iff } (t,s) < (r,q)$ lexicographically. We order
$\{\delta^{-1},\delta,\sigma_{ts},\ n\geq t>s\geq1\}^*$ by deg-lex
order.

Instead of $\sigma_{ij}$, we write simply $(i, j)$ or $( j, i)$. We
also set $(t_m, t_{m-1}, \dots, t_1)=\\(t_m, t_{m-1})(t_{m-1},
t_{m-2})\dots(t_2, t_1)$, where $t_j\neq t_{j+1},\ 1\leq j\leq m
-1$. In this notation, the defining relations of $B_n$ can be
written as
\begin{eqnarray*}
&&(t_3, t_2, t_1)=(t_2, t_1, t_3)=(t_1, t_3, t_2),\ t_3
>t_2>t_1,\\
&&(k, l)(i, j)=(i, j)(k, l),\ k
>l,\ i > j,\ k > i,
\end{eqnarray*}
where either $k > i > j >l$, or $k >l > i > j$.

Let us assume the following notation: $V_{[t_2,t_1]}$, where $n\geq
t_2>t_1\geq1$, is a positive word in $(k, l)$ such that $t_2\geq k >
l \geq t_1$. We can use any capital Latin letter with indices
instead of $V$, and any appropriate numbers (for example, $t_3, t_0$
such that $t_3
>t_0$) instead of $t_2, t_1$. We will also use the following notations:
$V_{[t_2-1,t_1]}(t_2, t_1)=(t_2, t_1)V'_{[t_2-1,t_1]},\ t_2
>t_1$, where $V'_{[t_2-1,t_1]}=(V_{[t_2-1,t_1]})|_{(k,l)\mapsto
(k,l), \ if \ l\neq t_1; \ (k,t_1)\mapsto(t_2,k)}$;
$W_{[t_2-1,t_1]}(t_1, t_0)=(t_1, t_0)W^\star_{[t_2-1,t_1]},\
t_2>t_1>t_0,$ where
$W^\star_{[t_2-1,t_1]}=(W_{[t_2-1,t_1]})|_{(k,l)\mapsto (k,l), \ if
\ l\neq t_1; \ (k,t_1)\mapsto(k,t_0)}$.

\begin{theorem}(\cite{Bo09})\
 A Gr\"{o}bner-Shirshov basis of the braid group $B_{n+1}$ in
the Birman-Ko-Lee generators consists of the following relations:
\begin{eqnarray*}
&&(k, l)(i, j) = (i, j)(k, l),\ k >l > i > j,\\
&&(k, l)V_{[j-1,1]}(i,j) = (i, j)(k, l)V_{[j-1,1]},\ k > i > j >l,\\
&&(t_3, t_2)(t_2, t_1)=(t_2, t_1)(t_3, t_1),\\
&&(t_3, t_1)V_{[t_2-1,1]}(t_3, t_2)=(t_2, t_1)(t_3, t_1)V_{[t_2-1,1]},\\
&&(t,s)V_{[t_2-1,1]}(t_2, t_1)W_{[t_3-1,t_1]}(t_3, t_1)=
(t_3, t_2)(t,s)V_{[t_2-1,1]}(t_2, t_1)W'_{[t_3-1,t_1]},\\
&&(t_3,s)V_{[t_2-1,1]}(t_2, t_1)W_{[t_3-1,t_1]}(t_3, t_1)=
(t_2,s)(t_3,s)V_{[t_2-1,1]}(t_2, t_1)W'_{[t_3-1,t_1]},\\
&&(2,1)V_{2[2,1]}(3,1)\dots V_{n-1[n-1,1]}(n,1)=\delta V'_{2[2,1]}\dots V'_{n-1[n-1,1]},\\
&&(t,s)\delta=\delta(t+1,s+1),\ \ (t,s)\delta^{-1}=\delta^{-1}(t-1,s-1),\ t\pm1,\ s\pm1\ (modn),\\
&&\delta\delta^{-1}=1,\ \delta^{-1}\delta=1,
\end{eqnarray*}
where $V_{[k,l]}$ means as before any word in $(i, j)$ such that
$k\geq i> j\geq l,\ t >t_3,\ t_2 > s$.
\end{theorem}

As results, we have the following corollaries.

\begin{corollary}
The semigroup of positive braids $BB^+_n$ in Birman-Ko-Lee
generators can be embedded into a group.
\end{corollary}

\begin{corollary}
The $S$-irreducible normal form of a word of $B_n$ in Birman-Ko-Lee
generators coincides with the Birman-Ko-Lee-Garside normal form
$\delta^kA,\ A\in BB^+_n$ of the word.
\end{corollary}

\subsection{Braid
groups in Adyan-Thurston generators}

The symmetry group is as follow:
\begin{eqnarray*}
S_{n+1}=gp\langle s_1,\dots, s_n\ |\ s_i^2=1, s_js_i=s_is_j\
(j-1>i), s_{i+1}s_is_{i+1}=s_is_{i+1}s_i\rangle.
\end{eqnarray*}
L. A. Bokut, L.-S.  Shiao \cite{bs} found the normal form for
$S_{n+1}$ in the following theorem.

\begin{theorem}(\cite{bs})
$N=\{s_{1i_1}s_{2i_2}\cdots s_{ni_n}|\ i_j\leq j+1\}$ is the
Gr\"{o}bner-Shirshov normal form for $S_{n+1}$ in generators
$s_i=(i, i+1)$ relative to the deg-lex order, where
$s_{ji}=s_js_{j-1}\cdots s_i\ (j\geq i),\ s_{jj+1}=1$.
\end{theorem}

Let $\alpha\in S_{n+1}$ and $
\overline{\alpha}=s_{1i_1}s_{2i_2}\cdots s_{ni_n}\in N$ be the
normal form of $\alpha$. Define the length of $\alpha$ as
$|\overline{\alpha}|=l(s_{1i_1}s_{2i_2}\cdots s_{ni_n})$ and
$\alpha\perp\beta$ if
$|\overline{\alpha\beta}|=|\overline{\alpha}|+|\overline{\beta}|$.
Moreover, each $\overline{\alpha}\in N$ has a unique expression
$\overline{\alpha}=s_{_{l_1i_{l_1}}}s_{_{l_2i_{l_2}}}\cdots
s_{_{l_ti_{l_t}}}$, where each $s_{_{l_ji_{l_j}}}\neq1$. Such a $t$
is called the breath of $\alpha$.

Now, we let
$$
B'_{n+1}=gp\langle r(\overline{\alpha}),\ \alpha\in
S_{n+1}\setminus\{1\}\ |\
r(\overline{\alpha})r(\overline{\beta})=r(\overline{\alpha\beta}),\
\alpha \perp \beta \rangle,
$$
where $r(\overline{\alpha})$ means a letter with the index
$\overline{\alpha}$.

Then for the braid group with $n+1$ generators (see \S3.1),
$B_{n+1}\cong B'_{n+1}$. Indeed, define $\theta: B_{n+1}\rightarrow
B'_{n+1},\ \sigma_i\mapsto r(s_i)$ and $\theta': B'_{n+1}\rightarrow
B_{n+1},\ r(\overline{\alpha})\mapsto \overline{\alpha}|_{s_i\mapsto
\sigma_i}$. Then two mappings are homomorphisms and
$\theta\theta'={l}_{B'_{n+1}}, \theta'\theta={l}_{B_{n+1}}$. Hence,
$$
B_{n+1}=gp\langle r(\overline{\alpha}),\ \alpha\in
S_{n+1}\setminus\{1\}\ |\
r(\overline{\alpha})r(\overline{\beta})=r(\overline{\alpha\beta}),\
\alpha \perp \beta \rangle.
$$
Let $X=\{r(\overline{\alpha}),\ \alpha\in S_{n+1}\setminus\{1\}\}$.
The generator $X$ of $B_{n+1}$ is called Adyan-Thurston generator.

Then the positive braid semigroup in generator $X$ is
$$
B_{n+1}^{+}=sgp\langle X\ |\
r(\overline{\alpha})r(\overline{\beta})=r(\overline{\alpha\beta}),\
\alpha \perp \beta \rangle.
$$

Let $s_1<s_2<\cdots<s_n$. Define
$r(\overline{\alpha})<r(\overline{\beta})$ if and only if
$|\overline{\alpha}|>|\overline{\beta}|$ or
$|\overline{\alpha}|=|\overline{\beta}|,\ \overline{\alpha}
<_{lex}\overline{\beta}$. It is clear that such an order on $X$ is
well-ordered. We will use the deg-lex order on $X^*$ in this
subsection.

\begin{theorem}(\cite{ChZhong-Braid})\label{t2.4}
A Gr\"{o}bner-Shirshov basis of $B_{n+1}^{+}$ in Adyan-Thurston
generator $X$ relative to the deg-lex order on $X^*$ is:
\begin{eqnarray*}
r(\overline{\alpha})r(\overline{\beta})&=&r(\overline{\alpha\beta}),\ \ \ \alpha \perp \beta, \\
r(\overline{\alpha})r(\overline{\beta\gamma})&=&r(\overline{\alpha\beta})r(\overline{\gamma}),\
\ \ \alpha \perp \beta \perp \gamma.
\end{eqnarray*}
\end{theorem}

\begin{theorem}(\cite{ChZhong-Braid})
A Gr\"{o}bner-Shirshov basis of $B_{n+1}$ in Adyan-Thurston
generator $X$ relative to the deg-lex order on $X^*$ is:
\begin{eqnarray*}
&&1)\ \ \ r(\overline{\alpha})r(\overline{\beta})=r(\overline{\alpha\beta}),\ \ \ \alpha \perp \beta, \\
&&2)\ \ \
r(\overline{\alpha})r(\overline{\beta\gamma})=r(\overline{\alpha\beta})r(\overline{\gamma}),\
\ \ \alpha \perp
\beta \perp \gamma,\\
&&3)\ \ \
r(\overline{\alpha})\Delta^{\varepsilon}=\Delta^{\varepsilon}r(\overline{\alpha}'),\
\ \
\overline{\alpha}'=\overline{\alpha}|_{s_i\mapsto s_{n+1-i}},\\
&&4)\ \ \ r(\overline{\alpha\beta})r(\overline{\gamma\mu})=\Delta
r(\overline{\alpha}')r(\overline{\mu}),\ \ \ \alpha \perp \beta
\perp \gamma \perp \mu,\ r(\overline{\beta\gamma})=\Delta,\\
&&5)\ \ \ \Delta^{\varepsilon}\Delta^{-\varepsilon}=1.
\end{eqnarray*}
\end{theorem}

\subsection{GSB and normal forms for free inverse semigroups}

We recall that an inverse semigroup is a semigroup in which every
element $a$ has a uniquely determined $a^{-1}$ such that
$aa^{-1}a=a\mbox{ and } a^{-1}aa^{-1}=a^{-1}$.  Let ${\cal {FI}}(X)$
be a free inverse semigroup generated by a set $X$,
$X^{-1}=\{x^{-1}|x\in X\}$ with $X\cap X^{-1}=\varnothing$. Denote
$X\cup X^{-1}$ by $Y$. Then ${\cal {FI}}(X)$ has the following
semigroup presentation
$$
\mathcal {FI}(X)=sgp\langle Y|\ aa^{-1}a=a,\
aa^{-1}bb^{-1}=bb^{-1}aa^{-1}, \ a,b\in Y^* \rangle
$$
where $1^{-1}=1,\ (x^{-1})^{-1}=x\ (x\in X)$ and $(y_1y_2\cdots
y_n)^{-1}=y_n^{-1}\cdots y_2^{-1}y_1^{-1}\ (y_1,\ y_2,\cdots,\
y_n\in Y)$.

Let us assume that the set $Y$ is well-ordered by an order $<$. Let
$<$ be also the corresponding deg-lex order of $Y^*$. For any
$u=y_1y_2\cdots y_n \ (y_1,\ y_2,\cdots,\ y_n\in Y)$, let
$fir(u)=y_1$.

We give inductively definitions in $Y^*$ of an idempotent, canonical
idempotent, prime canonical idempotent, ordered  (prime) canonical
idempotent and factors of a canonical idempotent, all of which but
(prime) idempotent and ordered (prime) canonical idempotent are
defined in O. Poliakova, B. M. Schein \cite{Schein}.

 (i) The empty word 1 is an idempotent, a canonical idempotent,
and an ordered canonical idempotent. This canonical idempotent has
no factors.

 (ii) If $h$ is an idempotent and $x\in Y$, then $x^{-1}hx$ is both an
idempotent and a prime idempotent. If $h$ is a canonical idempotent,
$x\in Y$ and the first letters of factors of $h$ are different from
$x$, then $x^{-1}hx$ is both a canonical idempotent and a prime
canonical idempotent. This canonical idempotent is its own factor.
Moreover, if the subword $h$ in this canonical idempotent is an
ordered canonical idempotent, then $x^{-1}hx$ is both an ordered
canonical idempotent and an ordered prime canonical idempotent.

 (iii) If $e_1, e_2,\cdots, e_m \ (m> 1)$ are prime idempotents,
then $e=e_1e_2\cdots e_m$ is an idempotent. Moreover, if $e_1,
e_2,\cdots, e_m$ are prime canonical idempotents and their first
letters are pairwise distinct, then $e=e_1e_2\cdots e_m$ is a
canonical idempotent and $e_1, e_2,\cdots, e_m$ are factors of $e$.
For this canonical idempotent, if $e_1, e_2,\cdots, e_m$ are ordered
canonical idempotents and $e\leq e_{i_1}e_{i_2}\cdots e_{i_m}$ for
any permutation $ (i_1, i_2,\cdots, i_m)$  of $ (1, 2,\cdots,m)$,
then $e$ is an ordered canonical idempotent.

\begin{theorem}(\cite{BoChZhao-inverse-sg})\label{gsb-sgp}
Let $X$ be a set, $X^{-1}=\{x^{-1}|x\in X\}$ with $X\cap
X^{-1}=\varnothing$, $<$ a well ordering on $Y=X\cup X^{-1}$ and
also the deg-lex order of $Y^*$, ${\cal {FI}}(X)$ the free inverse
semigroup generated by $X$. Let $S\subset k\langle Y\rangle$ be the
set of the following defining relations (a) and (b):

 (a) $ef=fe$, where both $e$ and $f$ are ordered prime canonical idempotents,
$ef$ is a canonical idempotent and $fe<ef$;

(b) $x^{-1}e'xf'x^{-1}=f'x^{-1}e'$, where $x\in Y$, both $x^{-1}e'x$
and $xf'x^{-1}$ are ordered prime canonical idempotents.

Then ${\cal {FI}}(X)=sgp\langle Y|S\rangle$ and $S$ is a
Gr\"{o}bner-Shirshov basis in $k\langle Y\rangle$.
\end{theorem}

By Theorem \ref{gsb-sgp} and CD-lemma for associative algebras, $Irr
(S)$ is normal forms of the free inverse semigroup $\mathcal
{FI}(X)$. It is easy to see that $Irr (S) = \{ u \in (X\cup
X^{-1})^* | u \neq a\bar{s}b ,s\in S,a ,b \in (X\cup X^{-1})^*\}$
consists of the word $u_0e_1u_1\cdots e_mu_m\in (X\cup X^{-1})^* $,
where $m\geq0, \ u_1,\cdots, u_{m-1}\neq1$, $u_0u_1\cdots u_m$ has
no subword of form $yy^{-1}$ for $y\in X\cup X^{-1}$, $e_1,\cdots,
e_m$ are ordered canonical idempotents, and the first (last,
respectively) letters of the factors of $e_i\ (1\leq i\leq m)$ are
not equal to the first (last, respectively) letter of $u_i\
(u_{i-1},\ \mbox{respectively})$. Thus $Irr (S)$ is a set of
canonical words in the sense of \cite{Schein}, and different words
in $Irr (S)$ represent different elements in $\mathcal {FI}(X)$.

\subsection{Embeddings of algebras}

 By using CD-lemma for associative algebras, in \cite{BoChMo-embed}, they give
another proofs for the following theorems (G. Higman, B. H. Neumann,
H. Neumann \cite{HNN49}, A. I. Malcev \cite{Ma}, T. Evans
\cite{Ev}): every countably generated group (resp. associative
algebra, semigroup) can be embedded into a two-generated group
(resp. associative algebra, semigroup). Also some new results are
proved.

\begin{theorem}(\cite{BoChMo-embed})\
(i) Every countably generated associative algebra over a countable
field $k$ can be embedded into a simple two-generated associative
algebra.

(ii) Every countably generated semigroup can be embedded into a
(0-)simple two-generated semigroup.
\end{theorem}

By using CD-lemma for Lie algebras over a field (see \S2.1), in
\cite{BoChMo-embed}, they give another proof of the Shirshov's
theorem (A. I. Shirshov \cite{Shir1}): every countably generated Lie
algebra can be embedded into a two-generated Lie algebra. Namely,
let $ L=Lie( X|S) $ is a Lie algebra generated by $X$ with relations
$S$, where $X=\{x_{i},i=1,2,\dots\}$ and $S$ is a GSB in the free
Lie algebra $Lie(X)$ on deg-lex order. Let $ H=Lie( X, a,b | S,
[aab^{i}ab]=x_{i}, \ i=1,2,\dots). $ Then $ \{ S, [aab^{i}ab]=x_{i},
\ i=1,2,\dots\} $ is a GSB in $Lie( X,a,b )$ on deg-lex ordering
with $a>b>x_i$. By CD-lemma for Lie algebras, $L$ is a subalgebra of
$H$ which is generated by $\{a,b\}$.

\begin{theorem}(\cite{BoChMo-embed})\
Every every countably generated Lie algebra over a countable field
$k$ can be embedded into a simple two-generated Lie algebra.
\end{theorem}

CD-lemma for associative differential algebras with unit is
established in  \cite{ChChLi-GSB-diff}, see \S2.6. By applying this
lemma, they prove the following theorem.

\begin{theorem}(\cite{BoChMo-embed})\
(i) Every countably generated associative differential algebra can
be embedded into a two-generated associative differential algebra.

(ii) Any associative differential algebra can be embedded into a
simple associative differential algebra.

(iii) Every countably generated associative differential algebra
with countable set $\mathcal{D}$ of differential operations over a
countable field $k$ can be embedded into a simple two-generated
associative differential algebra.
\end{theorem}

CD-lemma for associative algebra with multiple operations is
established in  \cite{BoChQiu-CD-Omega}, see \S2.2. By applying this
lemma, they prove the following theorems.

\begin{theorem}(\cite{BoChMo-embed})\
(i) Every countably generated associative $\Omega$-algebra can be
embedded into a two-generated associative $\Omega$-algebra.

(ii) Any associative $\Omega$-algebra can be embedded into a simple
associative $\Omega$-algebras.

(iii) Each countably generated associative $\Omega$-algebra with
countable multiple operations $\Omega$ over a countable field $k$
can be embedded into a simple two-generated associative
$\Omega$-algebra.
\end{theorem}

The definition of associative $\lambda$-differential algebra is
mentioned in \S2.2.

\begin{theorem}(\cite{BoChMo-embed})\
(i) Each countably generated associative $\lambda$-differential
algebra can be embedded into a two-generated associative
$\lambda$-differential algebra.

(ii) Each associative $\lambda$-differential algebra can be embedded
into a simple associative $\lambda$-differential algebra.

(iii) Each countably generated associative $\lambda$-differential
algebra over a countable field $k$ can be embedded into a simple
two-generated associative $\lambda$-differential algebra.
\end{theorem}

\subsection{Word problem for Novikov's
and Boone's group}

Let us call the set of letters $a_1,\ \cdots,a_n$ the principal
alphabet and refer to the letters
$$
q_1, \cdots,q_\lambda, r_1,\cdots, r_\lambda, l_1,\cdots,l_\lambda
$$
as signal ($ n ,\lambda>0$).\ Append a copy to the previous
alphabet, namely, the letters
$$
q_1^+,\cdots,\ q_\lambda^+,\ r_1^+,\
\cdots,r_\lambda^+,\l_1^+,\cdots,l_\lambda^+.
$$
Consider the set $\{(A_i,B_i)|1 \leq i\leq \lambda \}$ which is
constituted by pairs of nonempty words in the principal alphabet.
 Consider the chain of the
following groups: $G_0,\ G_1,\ G_2$ and $\mathcal {A}_{p_1p_2}$
given as follows.
\begin{eqnarray*}
G_0&=&gp\langle \ q_i,\ r_i,\ q_i^+,\ r_i^+\ (1\leq i\leq
\lambda) \ |  \emptyset  \rangle\\
G_1&=&gp\langle \ G_0,a_j,a_{j}^{+}\ (\ 1\leq j\leq n) \ |
\ q_ia_j=a_jq_iq_i,\\
&& \ \ \ \ r_ir_ia_j=a_jr_i,\ q_i^+q_i^+a_j^+=a_j^+q_i^+,\
r_i^+a_j^+=a_j^+r_i^+r_i^+ \ \rangle\\
G_2&=&gp\langle  \ G_1,\ l_i,\l_i^+\ (\ 1\leq i\leq \lambda \ ) \ |
\ a_jl_i=l_ia_j,\
a_j^+l_i^+=l_i^+a_j^+\ \ \rangle\\
\mathcal {A}_{p_1p_2}&=&gp\langle  \ G_2,\ p_1,\ p_2 \  |\
q_i^+l_i^+p_1l_iq_i=A_i^+p_1A_i,\
r_i^+p_1r_i=p_1,\\
&& \ \ \ \ r_il_ip_2l_i^+r_i^+=B_ip_2B_i^+,\ q_ip_2q_i^+=p_2 \
\rangle
\end{eqnarray*}

Then $\mathcal {A}_{p_1p_2}$ is called the Novikov group.

\ \

Now, the principal alphabet is constituted by letters ${s_{j}}^{,}s$
and ${q_{j}}^{,}s,\ \ 1\leq j\leq m$. The letters
$x,y,{l_{i}}^{,}s$, and ${r_{i}}^{,}s,\ \ 1\leq i\leq n$, form an
auxiliary alphabet. Denote by $\{(\Sigma_{i},\Gamma_{i}),\ \ 1\leq
i\leq n\}$ the set of pairs of special words, i.e., the words of the
form $sq_{j}s'$, where $s$ and $s'$ are words in the alphabet
$\{s_{j}\}$. Designate $q_{1}=q$. The Boone group $G(T,q)$ is given
by the above generators and the following defining relations:
\begin{description}
\item{(i)}\ \ \ $y^2s_{j}=s_{j}y,\  \ xs_{j}=s_{j}x
^2$;

\item{(ii)}\ \  \ $s_{j}l_{i}=yl_{i}ys_{j},\  \
s_{j}xr_{i}x=r_{i}s_{j}$;

\item{(iii)}\ \  \ $l_{i}\Gamma_{i}r_{i}=\Sigma_{i}$;

\item{(iv)}\ \  \ $l_{i}t=tl_{i},\  \ yt=ty$;

\item{(v)}\ \  \ $r_{i}k=kr_{i},\  \ xk=kx$;

\item{(vi)}\ \  \ $q^{-1}tqk=kq^{-1}tq$.

\end{description}

In the paper \cite{ChChLuo}, direct proofs of GSB for the Novikov
group and Boone group are given, respectively. Then the
corresponding normal forms are obtained. These proofs are different
from those in L. A. Bokut  \cite{bokut66,bokut67}.

\subsection{PBW basis of $U_q(A_N)$}

Let $A=(a_{ij})$ be an integral symmetrizable $N \times N$ Cartan
matrix so that $a_{ii}=2$, $a_{ij} \leq0 \ (i \neq j)$ and there
exists a diagonal matrix $D$ with diagonal entries $d_i$ which are
nonzero integers such that the product $DA$ is symmetric. Let q be a
nonzero element of $k$ such that $q^{4d_i} \neq 1$ for each $i$.
Then the quantum enveloping algebra is
$$
U_q(A)=k \langle X \cup H \cup Y|S^+ \cup K \cup T \cup S^- \rangle,
$$
where
\begin{eqnarray*}
X&=&\{x_i\},\\
H&=&\{h_i^{\pm1}\},\\
Y&=&\{y_i\},\\
S^+&=&\{ \sum_{\nu=0}^{1-a_{ij}}(-1)^\nu \left(
\begin{array}{c}
1-a_{ij} \\
\nu
\end{array}
 \right)_tx_i^{1-a_{ij}-\nu}x_jx_i^\nu, \ where \ i \neq j, \ t=q^{2d_i}\},\\
S^-&=& \{ \sum\limits_{\nu=0}^{1-a_{ij}}(-1)^\nu\left(
\begin{array}{c}
1-a_{ij} \\
\nu
\end{array}
\right)_ty_i^{1-a_{ij}-\nu}y_jy_i^\nu, \ where \ i \neq j, \ t=q^{2d_i}\},\\
K&=&\{h_ih_j-h_jh_i, h_ih_i^{-1}-1, h_i^{-1}h_i-1,
x_jh_i^{\pm1}-q^{\mp1}
d_ia_{ij}h^{\pm1}x_j, h_i^{\pm1}y_j-q^{\mp1}y_jh^{\pm1}\},  \\
T&=&\{x_iy_j-y_jx_i- \delta_{ij}
\frac{h_i^2-h_i^{-2}}{q^{2d_i}-q^{-2d_i}}\} \ \ \mbox{ and }\\
\left(
\begin{array}{c}
m \\
 n \\
 \end{array}
 \right)_t&=&\left\{
 \begin{array}{cc}
 \prod\limits_{i=1}^n \frac{t^{m-i+1}-t^{i-m-1}}{t^i-t^{-i}} & \ \ (for \ m>n>0),\\
 1 & \ \ \ (for \ n=0 \ or \ m=n).
 \end{array}
 \right.
\end{eqnarray*}

Let
$$
A=A_N=\left(
  \begin{array}{ccccc}
    2 & -1 & 0 & \cdots & 0 \\
    -1 & 2 & -1 & \cdots & 0 \\
    0 & -1 & 2 & \cdots & 0 \\
    \cdot & \cdot & \cdot & \cdot & \cdot \\
    0 & 0 & 0 & \cdots & 2 \\
  \end{array}
\right) \ \mbox{ and }\ \ q^8 \neq1.
$$

We introduce some new variables defined by Jimbo (see \cite{Ya})
which generate $U_q(A_N)$:
$$
\widetilde{X}=\{x_{ij}, 1 \leq i<j \leq N+1\},
$$
where
$$
x_{ij}=\left\{
\begin{array}{cc}
x_i & \ \ \  j=i+1,\\
qx_{i,j-1}x_{j-1,j}-q^{-1}x_{j-1,j}x_{i,j-1} & \ \ \ j>i+1.
\end{array}
\right.
$$
We now order the set $\widetilde{X}$ in the following way.
$$
x_{mn}>x_{ij} \Longleftrightarrow(m,n)>_{lex}(i,j).
$$
Let us recall from Yamane \cite{Ya}  the following notation:
\begin{eqnarray*}
C_1&=&\{((i,j),(m,n))|i=m<j<n\},\\
C_2&=&\{((i,j),(m,n))|i<m<n<j\},\\
C_3&=&\{((i,j),(m,n))|i<m<j=n\},\\
C_4&=&\{((i,j),(m,n))|i<m<j<n\},\\
C_5&=&\{((i,j),(m,n))|i<j=m<n\},\\
C_6&=&\{((i,j),(m,n))|i<j<m<n\}.
\end{eqnarray*}
Let the set $\widetilde{S}^+$ consist of Jimbo relations:
\begin{eqnarray*}
x_{mn}x_{ij}&-&q^{-2}x_{ij}x_{mn}  \ \ \ \ \ \ \ \ \ \ \ \ \ \ \ \
\ \ \ \ \ \ \ \ \ ((i,j),(m,n)) \in C_1 \cup C_3,\\
x_{mn}x_{ij}&-&x_{ij}x_{mn} \ \ \ \ \ \ \ \ \ \ \ \ \ \ \ \ \ \ \ \
\ \ \ \ \ \ \ \ \ ((i,j),(m,n)) \in C_2
\cup C_6,\\
x_{mn}x_{ij}&-&x_{ij}x_{mn}+(q^2-q^{-2})x_{in}x_{mj} \ \ \
((i,j),(m,n)) \in C_4,\\
x_{mn}x_{ij}&-&q^2x_{ij}x_{mn}+qx_{in} \ \ \ \ \ \ \ \ \ \ \ \ \ \ \
\ \ ((i,j),(m,n)) \in C_5.
\end{eqnarray*}
It is easily seen that $U^+_q(A_N)=k\langle
\widetilde{X}|\widetilde{S^+}\rangle$.

\ \

In the paper \cite{ChShaoShum}, a direct proof is given that
$\widetilde{S}^+$ is a Gr\"{o}bner-Shirshov basis for $k\langle
\widetilde{X}|\widetilde{S^+}\rangle=U^+_q(A_N)$ (\cite{BoMa}). The
proof is different from one in L. A. Bokut, P. Malcolmson
\cite{BoMa}.

\subsection{GSB for free dendriform algebras}

The $L$-algebra and the dendriform algebra are mentioned in \S2.5
and \S2.3, respectively.

Let $L(X)$ be the free $L$-algebra generated by a set $X$ and $N$
the set of all $L$-words, see \S2.5. Let the order on $N$ be as in
\S2.5.

It is clear that the free dendriform algebra generated by $X$,
denoted by $DD(X)$, has an expression
\begin{eqnarray*}
L(X&|&(x \prec y)\prec z = x \prec (y \prec z) + x \prec (y \succ
z),\\
&&  x \succ(y \succ z) = (x \succ y)\succ z + (x \prec y)\succ z,\
x,y,z\in N).
\end{eqnarray*}

The following theorem gives a GSB for $DD(X)$.

\begin{theorem} (\cite{ChWang-GSB-Den})
 Let the order on
$N$ be as before. Let
\begin{eqnarray*}
f_1(x,y,z)& =&(x\prec y)\prec z - x \prec (y \prec z) -x \prec (y
\succ z),\\
f_2(x,y,z)&=&(x\prec y)\succ z + (x \succ y)\succ z-x\succ (y\succ
z),\\
f_3(x,y,z,v)&=&((x\succ y)\succ z)\succ v - (x\succ y)\succ (z\succ
v) + (x\succ(y\prec z))\succ v.
\end{eqnarray*}
Then, $S=\{f_1(x,y,z),~f_2(x,y,z),~f_3(x,y,z,v)|~x,y,z,v \in N \}$
is a Gr\"{o}bner-Shirshov basis in $L(X)$.
\end{theorem}

Then, by using the CD-lemma for $L$-algebras, see
\cite{ChHuang-L-alg}, $Irr(S)$ is normal words of free dendriform
algebra $DD(X)$.

\subsection{Anti-commutative GSB of a free Lie algebra
relative to Lyndon-Shirshov words}

Let $X$ be a well-ordered set,  $>_{lex}$, $>_{dex-lex}$
lexicographical order and degree-lexicographical order on $X^{*}$
respectively.

Let $>_{n-lex}$ be the nonassociative lexicographical order on
$X^{**}$. For example, let $X=\{x_i|1\leq i\leq n \mbox{ and }
x_1<x_2<\ldots<x_n\}$. Then $x_1(x_2x_3)>_{n-lex}(x_1x_2)x_3$ since
$x_1>_{n-lex}(x_1x_2)$.

Let $>_{n-dex-lex}$ be the  nonassociative degree-lexicographical
order on $X^{**}$. For example,
$x_1(x_2x_3)<_{n-deg-lex}(x_1x_2)x_3$ since
$x_1<_{n-deg-lex}(x_1x_2)$ (the degree of $x_1$ is one and the
degree of $(x_2x_3)$ is two).

An order $\succ_1$ on $X^{**}$:

$$
(u)\succ_1(v)\Leftrightarrow either \ \ \ u>_{lex}v\ \ \  or\ \ \
(u=v\ \ and \ \ (u)>_{n-dex-lex}(v)).
$$

An order $\succ_2$ on $X^{**}$:

$$
(u)\succ_2(v)\Leftrightarrow either \ \ \ u>_{deg-lex}v\ \ \  or\ \
\ (u=v\ \ and \ \ (u)>_{n-dex-lex}(v)).
$$

We define normal words $(u)\in X^{**}$ by induction:\\

(i) $x_i\in X$ is a normal word. \\

(ii) $(u)=(u_1)(u_2)$ is normal if and only if both $(u_1)$ and
$(u_2)$ are
normal and $(u_1)\succ_1(u_2)$.\\

Denote $N$ the set of all normal words $[u]$ on $X$.

Let  $AC(X)$ be a $k$-space spanned by $N$. Now define the product
of normal words by the following way: for any
 $[u],[v]\in N$,
\begin{equation*}
\lbrack u][v]=\left\{
\begin{array}{r@{\quad}l}
[[u][v]], & \ \ \ \mbox{ if } \  [u]\succ_1[v] \\
-[[v][u]], & \ \ \ \mbox{ if } \   [u]\prec_1[v] \\
0,\text{ \ \ } & \ \ \ \mbox{ if } \   [u]=[v]%
\end{array}%
\right.
\end{equation*}

\noindent \textbf{Remark} \ By definition, for any $(u)\in X^{**}$,
there exists a unique $[v]\in N$ such that, in $AC(X)$, $(u)=\pm
[v]$ or 0.

Then  $AC(X)$ is a free anti-commutative algebra generated by $X$.

The order $\succ_2$ is a monomial  order on $N$.

Shirshov \cite{Sh62a} proved a CD-lemma for free anti-commutative
algebra $AC(X)$, where there is only one composition: inclusion. By
using  the above order $\succ_2$ on $N$, Shirshov's (\cite{Sh62a})
CD-lemma holds.

\begin{theorem}(\cite{BoChLi10}) Let the order $\succ_2$ on $N$ be defined as before and
\begin{center}
$S=\{([[u]][[v]])[[w]]-([[u]][[w]])[[v]]-[[u]]([[v]][[w]]) \ |
~u>_{lex}v>_{lex}w\ and \ [[u]],[[v]],[[w]]\ are\ NLSWs\}.$
\end{center}
 Then
$S$ is a Gr\"obner-Shirshov basis in $AC(X)$.
\end{theorem}

By CD-lemma for anti-commutative algebras, $Irr(S)=NLSW(X)$ is
linear basis of the free Lie algebra $Lie(X)$ which is due to A. I.
Shirshov \cite{Shir1} and K.-T. Chen, R. Fox, R. Lyndon \cite{CFL}.

\subsection{Free
partially commutative groups, associative algebras and Lie algebras}

Let $X$ be a set, $Free(X)$ the free algebra over a field $k$
generated by $X$, for example, free monoid, group, associative
algebra over a field, Lie algebra and so on. Let
$\vartheta\subseteq(X\times X)\backslash \{(x,x) \mid x\in X\}$.
Then $Free(X|R)=Free(X)/Id(R)$ with generators $X$ and defining
relations $R=\{ab=ba,\ (a,b)\in\vartheta\}$ (in the Lie case
$R=\{(a,b)=0,\ (a,b)\in\vartheta\})$ is the free partially
commutative algebra.

Denote by
\begin{enumerate}
\item[---] $gp\langle X|
\vartheta\rangle$ the free partially commutative group.
\item[---] $sgp\langle X|
\vartheta\rangle$ the free partially commutative monid.
\item[---] $k\langle X|
\vartheta\rangle$ the free partially commutative associative
algebra.
\item[---] $Lie( X|\vartheta)$ the free partially commutative Lie
algebra.
\end{enumerate}

Let  $<$ be a well ordering on $X$.  Throughout this subsection,  if
$a>b$ and $(a,b)\in \vartheta$ or $(b,a)\in \vartheta$, we denote
$a\rhd b$. Generally, for any set $Y$, $a\rhd Y$ means $a\rhd y$ for
any $y\in Y$. For any $u=x_{i_1}\cdots x_{i_n}\in X^*$ where
$x_{i_j}\in X$, we denote the set $\{x_{i_j}, j=1,\ldots, n\}$ by
$supp(u)$.

For free partially commutative group $gp\langle X|
\vartheta\rangle$, let $\prec$ be a well ordering on $X$.  We extend
this order to $X\cup X^{-1}$ as follows: for any $x, y\in X,\
\varepsilon, \eta=\pm1$
\begin{enumerate}
\item [(i)] $x>x^{-1}$, if $x\in X$;
\item [(ii)] $x^{\varepsilon}>y^{\eta}\Longleftrightarrow x\succ
y$.
\end{enumerate}

For any $a, b\in X$, if $a^{\varepsilon}>b^{\eta}\ (\varepsilon,
\eta=\pm1)$ and $(a,b)\in \vartheta$  or $(b,a)\in \vartheta$, we
denote $a^{\varepsilon}\rhd b^{\eta}$. Moreover, for any set $Y$,
$a\rhd Y$ means $a\rhd y$ for any $y\in Y$.

It is obvious that
$$
gp\langle X| \vartheta\rangle=sgp\langle X\cup
X^{-1}|x^{\varepsilon}x^{-\varepsilon}=1, \ a^{\varepsilon}
b^{\eta}=b^{\eta}a^{\varepsilon},\ x\in X,\  (a,b)\in \vartheta,\
\varepsilon, \eta=\pm1\rangle.
$$

By using CD-lemma for associative algebras, the following results
are proved in \cite{ChMo-Partialcomm}.

\begin{theorem}(\cite{bs})\
With the deg-lex order on $X^*$, the set $S=\{xuy-yxu \mid x,y\in X,
u\in X^*, x\rhd y\rhd supp(u)\}$ forms a Gr\"{o}bner-Shirshov basis
of the free partially commutative associative algebra $k\langle X|
\vartheta\rangle$. As a result, $Irr(S) = \{ u \in X^* |  u \neq
a\bar{s}b,s\in S,a ,b \in X^*\}$ is a $k$-basis of $k\langle X|
\vartheta\rangle$. Then,  $Irr(S)$ is also a normal form of the free
partially commutative monoid $sgp\langle X| \vartheta\rangle$.
\end{theorem}

\begin{theorem}(\cite{Esyp})\
Let the notations be as the above. Then with the deg-lex order on
$(X\cup X^{-1})^*$, the set $S=\{x^{\varepsilon}uy^{\eta}=
y^{\eta}x^{\varepsilon}u,\ z^{\gamma}z^{-\gamma}=1\mid \varepsilon,
\eta, \gamma=\pm1, \ x,y,z\in X\cup X^{-1},\  u\in (X\cup X^{-1})^*,
\ x^{\varepsilon}\rhd y^{\eta}\rhd supp(u)\}$ forms a
Gr\"{o}bner-Shirshov basis of the free partially commutative group
$gp\langle X| \vartheta\rangle$. As a result, $Irr(S) = \{ u \in
(X\cup X^{-1})^* | u \neq a\bar{s}b,s\in S,a ,b \in (X\cup
X^{-1})^*\}$ is a normal form of $gp\langle X| \vartheta\rangle$.
\end{theorem}

By using CD-lemma for Lie algebras, we have

\begin{theorem}\label{t3.9.3}(\cite{ChMo-Partialcomm})\
With deg-lex order on $X^*$, the set $S=\{[xuy] \ | \ x,y,z\in X,\
u\in X^*, x\rhd y\rhd supp(u)\}$ forms a Gr\"{o}bner-Shirshov basis
of the partially commutative Lie algebra $Lie(X| \vartheta)$. As a
result, $Irr(S)=\{[u]\in NLSW(X) \ | \ u\neq{a\bar{s}b}, \ s\in{S},\
a,b\in{X^*}\}$ is a $k$-basis of $Lie(X| \vartheta)$.
\end{theorem}

Theorem \ref{t3.9.3} is also proved by E. Poroshenko
\cite{poroshenko}.

\subsection{Plactic monoids in row generators}

Let $X=\{x_1,\dots,x_n\}$ be a set of $n$ elements with the order
$x_1<\dots<x_n$. The plactic monoid (see M. Lothaire
\cite{Lothaire}, Chapter 5) on the alphabet $X$ is $P_n=sgp\langle
X|T\rangle$, where $T$ consists of the Knuth relations
\begin{eqnarray*}&&x_ix_kx_j= x_kx_ix_j\ \ \ (x_i\leq x_j<x_k),\\
&&x_jx_ix_k= x_jx_kx_i\ \ \ (x_i< x_j\leq x_k).
\end{eqnarray*}

A nondecreasing word $R\in X^*$ is called a row, for example,
$x_1x_1x_3x_5x_5x_5x_6$ is a row when $n\geq6$. Let $Y$ be the set
of all rows in $X^*$. For any $R=r_1\dots r_u,\ S=s_1\dots s_v\in
Y,\ r_i,s_j\in X$, we say $R$ dominates $S$ if $u\leq v$ and for
$i=1, \dots, u,\ r_i>s_i$.

The multiplication of two rows is defined by Schensted's algorithm:
 for a row $R$ and $x\in X$,
\[
R\cdot x= \left\{
\begin{array}{cc}
    Rx,& if~Rx~is~a~row \\
    y\cdot R',& otherwise
\end{array}
\right.
\]
where $y$ is the leftmost letter in $R$ which is strictly larger
than $x$, and $R'=R|_{y\mapsto x}$.

Then, for any $R,S\in Y$, it is clear that there exist uniquely
$R',S'\in Y$ such that $R\cdot S=R'\cdot S'$ in $P_n$, where $R'$
dominates $ S'$, i.e., $|R'|\leq|S'|$ and each letter of $R'$ is
larger than the corresponding letter of $S'$.

We express the $P_n$ as follows:
$$
P_n=sgp\langle Y|R\cdot S=R'\cdot S',\ R,S\in Y\rangle.
$$

For any $R,S\in Y$, we define $R>S$ by deg-lex order on $X^*$. Then
we have a well ordering on $Y$.

It follows from CD-lemma for associative algebras that with the
deg-lex order on $Y^*$, the set $\{R\cdot S=R'\cdot S'| R,S\in Y\}$
is a GSB in $k\langle Y\rangle$, see M. Lothaire \cite{Lothaire},
Chapter 5. For $n=4$, we give a direct proof. Namely, any
composition $(f,g)_w=(RS)T-R(ST)$ is trivial, where $f=RS-R'S',\
g=ST-S^{''}T',\ w=RST$.

Let $R=1^{r_1}2^{r_2}3^{r_3}4^{r_4}$ be a row in $P_4$, where each
$r_i$ is non-negative integer. For convenience, denote by
$R=(r_1,r_2,r_3,r_4)$. Then for three rows $R=(r_1,r_2,r_3,r_4),\
S=(s_1,s_2,s_3,s_4),\ T=(t_1,t_2,t_3,t_4)$, we have

\begin{eqnarray*}
R(ST)&=&\left( \begin{array}{c}
               r_1~~~r_2~~~r_3~~~r_4 \\
               \left( \begin{array}{cccc}
               s_1 & s_2 & s_3 & s_4 \\
               t_1 & t_2 & t_3 & t_4
               \end{array} \right)
               \end{array} \right)
        =\left( \begin{array}{c}
                \left( \begin{array}{cccc}
                r_1 & r_2 & r_3 & r_4 \\
                0   & a_2 & a_3 & a_4
                \end{array} \right) \\
                b_1~~~b_2~~~b_3~~~b_4
                \end{array} \right)\\
    &=&\left( \begin{array}{c}
        0~~~0~~~c_3~~~c_4 \\
        \left( \begin{array}{cccc}
        d_1 & d_2 & d_3 & d_4 \\
        b_1 & b_2 & b_3 & b_4
        \end{array} \right)
        \end{array} \right)
    =\left( \begin{array}{cccc}
    0  & 0  & c_3 & c_4 \\
    0  & e_2 & e_3 & e_4 \\
    f_1 & f_2 & f_3 & f_4
    \end{array} \right)\\
(RS)T&=&\left( \begin{array}{c}
               \left( \begin{array}{cccc}
               r_1 & r_2 & r_3 & r_4 \\
               s_1 & s_2 & s_3 & s_4
               \end{array} \right) \\
               t_1~~~t_2~~~t_3~~~t_4
               \end{array} \right)
        =\left( \begin{array}{c}
                0~~~g_2~~~g_3~~~g_4 \\
                \left( \begin{array}{cccc}
                h_1 & h_2 & h_3 & h_4 \\
                t_1 & t_2 & t_3 & t_4
                \end{array} \right)\\
                \end{array} \right)\\
    &=&\left( \begin{array}{c}
        \left( \begin{array}{cccc}
        0  & g_2 & g_3 & g_4 \\
        0  & i_2 & i_3 & i_4
        \end{array} \right)\\
        j_1~~~j_2~~~j_3~~~j_4
        \end{array} \right)
    =\left( \begin{array}{cccc}
    0  & 0  & k_3 & k_4 \\
    0  & l_2 & l_3 & l_4 \\
    j_1 & j_2 & j_3 & j_4
    \end{array} \right)
\end{eqnarray*}

Now, one needs only to prove $l_2=e_2$, $l_3=e_3$, $l_4=e_4$,
$c_3=k_3$ and $c_4=k_4$. Let us mention that it takes almost 20
pages to calculate.

As the result, in \cite{ChLiJing-plactic}, a new direct elementary
proof (but not simple) of the following theorem is given.

\begin{theorem}
Let the notations be as above. Then with deg-lex order on $Y^*$,
$\{R\cdot S=R'\cdot S'| R,S\in Y\}$ is a Gr\"{o}bner-Shirshov basis
in $k\langle Y\rangle$ if $n=4$.
\end{theorem}

\subsection{Plactic algebras in standard generators }

Let  $P_n=sgp\langle X|T\rangle$ be the plactic monoid with
$n$-generators, see \S3.12. Then the semigroup algebra $k\langle
X|T\rangle$ is called the plactic algebra. For any $x_{i_1}\cdots
x_{i_t}\in X^*$, we denote $x_{i_1}\cdots x_{i_t}$ by
${i_1}\cdots{i_t}$.

In the paper \cite{Okninskii}, the following theorem is given.

\begin{theorem}(\cite{Okninskii}, Theorem 1, Theorem 3)
Let the notation be as above. Then with the deg-lex order on $X^*$,
the following statements hold.
\begin{enumerate}
\item[(i)] If $|X|=3$, then the set $S=\{332-323,\ 322-232,\ 331-313,\
311-131,\ 221-212,\ 211-121,\ 231-213,\ 312-132,\ 3212-2321,\
32131-31321,\ 32321-32132\}$ is a Gr\"{o}bner-Shirshov basis in
$k\langle x_1,x_2,x_3\rangle$.
\item[(ii)] If $|X|>3$, then the Gr\"{o}bner-Shirshov complement $T^c$ of
$T$ is infinite.
\end{enumerate}
\end{theorem}

 \subsection{Filtrations and distortion in
infinite-dimensional algebras}

Let $A$ be a linear algebra over a field $k$. An ascending
filtration $\alpha=\{A_n\}$ on $A$ is a sequence of subspaces
$A_0\subset A_1\subset\ldots\subset A_n\subset\ldots$  such that
$A=\sum_{n=0}^{\infty}A_n$ and $A_kA_l\subset A_{k+l}$ for all $k, l
= 0, 1, 2,\ldots$.  Given $a\in A$, the $\alpha$-degree of $a$,
denoted by $deg_{\alpha}a$, is defined as the least $n$ such that
$a\in A_n$. If $B$ is a subalgebra of $A$ with a filtration
$\alpha$, as above, then we call the filtration $\beta=\{B\cap
A_n\}$ the restriction of $\alpha$ to $B$ and we write $\beta=\alpha
\cap B$. In the case of monoids the terms of filtrations are simply
subsets, with all other conditions being the same. In the case of
groups the terms of filtrations must be closed under inverses.

A generic example is as follows. Let $A$ be a unital associative
algebra generated by a finite set $X$. Then a filtration
$\alpha=\{A_n\}$ arises on $A$, if one sets $A_0 = Span \{1\}$ and
$A_n = A_{n-1} + Span \{X^n\}$, for any $n>1$. The $\alpha$-degree
of $a\in A$ in this case is an ¡°ordinary¡± degree with respect to
the generating set $X$, that is the least degree of a polynomial in
$X$ equal $a$. We write $deg_{\alpha}a= deg_X a$. Such filtration
$\alpha$ is called the degree filtration defined by the generating
set $X$.

\begin{definition}(\cite{Olshanskii})\ Given two filtrations
$\beta=\{B_n\}$ and $\beta'=\{B_n'\}$ on the same
algebra $B$, we say that $\beta$ is majorated by $\beta'$ if there
is an integer $t>0$ such that $B_n\subseteq B_{tn}'$, for all
$n\geq0$. We then write $\beta\preceq\beta'$. If
$\beta\preceq\beta'$ and $\beta'\preceq\beta$ then we say that
$\beta$ and $\beta'$ are equivalent and write $\beta\sim\beta'$.

If $B$ is a finitely generated subalgebra of a finitely generated
algebra $A$ then we say then $B$ is embedded in $A$ without
distortion (or that $B$ is an undistorted subalgebra of $A$) if a
degree filtration of $B$ is equivalent to the restriction to $B$ of
a degree filtration of $A$
\end{definition}

\begin{definition}(\cite{Olshanskii})\
A filtration $\alpha=\{A_n\}$ on an algebra $A$  is called a tame
filtration if it satisfies: there is $c>0$ such that $dim(A_n) <
c^n$, for all $n = 1, 2, \dots$.
\end{definition}

In the paper \cite{Olshanskii}, the following questions are
discussed.

(1) Is it true that every tame filtration of an algebra $B$ is
equivalent (or equal) to a filtration restricted from the degree
filtration of a finitely generated algebra $A$ where $B$ is embedded
as a subalgebra?

(2) If the answer to the previous question is ``yes", can one choose
$A$ finitely presented? If not, indicate conditions ensuring that
the answer is still ``yes".

By using CD-lemmas for associative and Lie algebras, Y. Bahturin, A.
Olshanskii \cite{Olshanskii} proved the following theorems.

\begin{theorem}(\cite{Olshanskii}, Theorem 7)
Let $B$ be a unital associative algebra over a field $F$.

(1) A filtration $\beta$ on $B$ is tame if and only if
$\beta\sim\alpha\cap B$ where $\alpha$ is a degree filtration on a
unital 2-generator associative algebra $A$ where $B$ is embedded as
a subalgebra.

(2) A filtration $\beta$ on $B$ is tame if and only if
$\beta=\alpha\cap B$ where $\alpha$  is a degree filtration on a
unital finitely generated associative algebra $A$ where $B$ is
embedded as a subalgebra.
\end{theorem}

\begin{theorem}(\cite{Olshanskii}, Theorem 8)
Let $N$ be a countable monoid.

(i) There exist a 3-generator monoid $M$ where $N$ is embedded as a
submonoid.

(ii) If $N$ is finitely generated then the embedding of $N$ in a
3-generator monoid $M$ can be done without distortion.

(iii) A filtration $\beta$ on $N$ is a tame filtration if and only
if there is a finitely generated monoid $M$ with a degree filtration
$\alpha$ such that $\beta\sim\alpha\cap N$.
\end{theorem}

\begin{theorem}(\cite{Olshanskii}, Theorem 9)\
(1) A filtration $\chi$ on a Lie algebra $H$ is tame if and only if
$\chi\sim\gamma\cap H$ where $\chi$
 is the degree filtration on a 2-generator Lie algebra
$G$ where $H$ is embedded as a subalgebra, if and only if
$\chi\sim\rho\cap H$ where $\rho$ is the degree filtration on a
2-generator associative algebra $R$ where $H$ is embedded as a Lie
subalgebra.

(2) A filtration  $\chi$ on a Lie algebra $H$ is tame if and only if
$\chi=\gamma\cap H$ where $\gamma$ is the degree filtration on a
finitely generated Lie algebra $G$ where $H$ is embedded as a
subalgebra, if and only if $\chi=\rho\cap H$  where $\rho$ is the
degree filtration on a finitely generated associative algebra  $R$
where $H$ is embedded as a Lie subalgebra.
\end{theorem}

\begin{theorem}(\cite{Olshanskii}, Theorem 10)
Any finitely generated associative, respectively, Lie algebra can be
embedded without distortion in a simple 2-generator associative,
respectively, Lie algebra.
\end{theorem}

\begin{theorem}(\cite{Olshanskii}, Theorem 11)
Any finitely generated monoid $M$ with a recursively enumerable set
of defining relations  can be embedded in a finitely presented
monoid $M'$ as an undistorted submonoid.
\end{theorem}

\begin{theorem}(\cite{Olshanskii}, Theorem 13)
Let $B$ be an arbitrary  finitely generated unital associative
algebra with a recursively enumerable set of defining relations,
over a field $k$ which is finitely generated over prime subfield.
Then there exists a  finitely presented unital associative
$k$-algebra $A$ in which $B$ is contained as an undistorted unital
subalgebra.
\end{theorem}

\subsection{Sufficiency
conditions for Bokut' normal forms}

L. A. Bokut \cite{bokut66,bokut67} (see also L. A. Bokut, D. J.
Collins \cite{BoCollins}) introduced a method for producing normal
forms, called Bokut's normal forms, for groups obtained as a
sequence of HNN extensions starting with a free group. K. Kalorkoti
\cite{Kalorkoti82,Kalorkoti06,Kalorkoti09} has used these in several
applications. In any application we need to prove first of all that
normal forms exists, i.e., a certain rewriting process terminates.
Uniqueness is then established separately, it is not guaranteed.
Termination is usually ensured for both requirements that lead to a
fairly uniform approach.

In the paper \cite{Kalorkoti10}, the author provides sufficient
conditions for the existence and uniqueness of normal forms of
sequences of HNN extentions defined by Bokut'. Furthermore, the
author shows that under an assumption, which holds for various
applications, such normal forms always exist (but might not be
unique). The conditions are amenable to be used in automatic theorem
provers. It is discussed also how to obtain a Gr\"{o}bner-Shirshov
basis from the rewrite rules of Bokut's normal forms under certain
assumptions. An application drawn from a paper of S. Aaanderaa, D.
E. Cohen \cite{Aaanderaa} to illustrate to sufficiency conditions is
also given.

\subsection{Quantum groups of type $G_2$ and $D_4$}

 First, we recall some basic notions about Ringel-Hall algebra from
 \cite{[R1]}.

 Let $k$ be a finite field with $q$ elements, $\Lambda$ an
$k-$algebra. By
 $\Lambda-$mod we denote the category of finite dimensional right
 $\Lambda-$modules. For $M,N_1,\cdots,N_t\in \Lambda-$mod, let
 $F^M_{N_1,\dots,N_t}$ be the number of filtrations
 $$M=M_0\supseteq M_1\supseteq \cdots \supseteq M_{t-1}\supseteq M_t=0,$$
 such that $M_{i-1}/M_{i}\cong N_i$ for all $1\leq i\leq t.$

 Now, let $\Lambda$ be a finitary $k-$algebra, i.e., for any $M,N\in
 \Lambda-$mod, $\hbox{Ext}^1_{\Lambda}(M,N)$ is finite dimensional $k-$vector
 space. For each $M\in \Lambda-$mod, we denote by $[M]$ the isomorphism
 class of $M$ and by ${\mathbf{dim}}M$ the dimension vector of the $\Lambda$-module $M.$
 We have the well-known Euler form $\langle-,-\rangle$ defined by
 $$\langle {\mathbf{dim}}M,{\mathbf{dim}}N\rangle=
  \sum_{i=0}^{\infty}\mathbf{dim}\hbox{Ext}^i_\Lambda(M,N).$$
  Note that $(-,-)$ is the symmetrization of $\langle -,- \rangle.$
  It is well-known  that if $\Lambda$ is Dynkin type, that is one of the types
  $A_n, B_n, C_n, D_n, E_6, E_7, E_8,  F_4$ and $G_2,$
  then  $\langle {\mathbf{dim}}M,{\mathbf{dim}}N\rangle=
  {\mathbf{dim}}\hbox{Hom}_\Lambda(M,N)-{\mathbf{dim}}\hbox{Ext}^1_\Lambda(M,N).$

  The (twisted) Ringel-Hall algebra $\mathcal{H}(\Lambda)$ is the free $\mathbb{Q}(\upsilon)$-algebra
  with bases $\{u_{[M]}\mid M\in \Lambda-\hbox{mod}\}$
  indexed by the set of isomorphism classes of $\Lambda$-modules,
 with multiplication given by
 $$
 u_{[M]}u_{[N]}=\upsilon^{\langle \mathbf{dim}M,\mathbf{N}\rangle}
 \sum_{[L]}F^L_{M,N}u_{[L]}\hspace{1cm} \hbox{for\; all } M,N\in \Lambda-\hbox{mod},
 $$
 where $\mathbb{Q}$ is the field of quotient numbers, $\upsilon$ an
 indeterminate and $\mathbb{Q}(\upsilon)$ the factor field of the
 polynomial algebra $\mathbb{Q}[\upsilon]$.

   Then $\mathcal{H}(\Lambda)$ is an associative algebra with identity
  $1=u_0$, where $0$ denotes the isomorphism class of the trivial
  $\Lambda-$module 0.

  Let $\Delta$ be the graph defined by the symmetrizable Cartan matrix
  A and $\overrightarrow{\Delta}$ be obtained by choosing some
  orientation of the edges of the graph $\Delta$. Choose a
  $k-$species $\mathfrak{S}=(F_i,_iM_j)$ of type
  $\overrightarrow{\Delta}$ (see \cite{[R2]}). Then the  main result in \cite{[R1]} is

\begin{theorem}\label{t3.14.1}(\cite{[R1]})\
The map  $\eta : U^+_q(A) \longrightarrow
\mathcal{H}(\overrightarrow{\Delta})$  given by
  $$\eta(E_i)=[S_i]$$
  is an $\mathbb{Q}(\upsilon)-$algebra isomorphism, where $[S_i]$ is the isomorphism class of i'th simple
  representations of $\mathfrak{S}$-module.
\end{theorem}

  We take the set $B=\{u_{[M]}\mid M \ \hbox{is}\ \hbox{
  indecomposable}\}$ as a generating set for the Ringel-Hall algebra
  $\mathcal{H}(\overrightarrow{\Delta})$
  and denote by $T^+$ the set of all skew-commutator relations between the elements of
  $B.$ Then by the isomorphism $\eta $ in Theorem \ref{t3.14.1},
  we get a corresponding set $S^{+c}$ of relations in the positive
  part of quantum group $ U^+_q(A)$.

  Let  $\Delta=G_2\ (\hbox{resp.}\ D_4)$ and $\overrightarrow{\Delta}$ be obtained
  by choosing the following orientations:
for $G_2:$
       
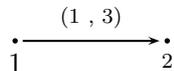
\begin{center}
\begin{pspicture}(0,0)(4,1)
 \psdot*[dotsize=2pt](1,1)
\psdot*[dotsize=2pt](3,1)
 \psline{->}(1.1,1)(2.9,1)
  \uput[d](1,1){\small {1}}
 \uput[d](3,1){\scriptsize {2}}
  \uput[u](2,1){\scriptsize {(1 , 3)}}
\end{pspicture}
\end{center}
and for $D_4:$
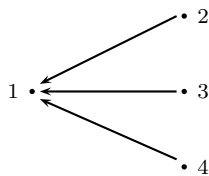
\begin{center}
\begin{pspicture}(0,-0.5)(4,2)
 \psdot*[dotsize=2pt](1,1)
\psdot*[dotsize=2pt](3,1)
 \psdot*[dotsize=2pt](3,2)
 \psdot*[dotsize=2pt](3,0)
 \uput[l](1,1){\scriptsize {1}}
 \uput[r](3,1){\scriptsize {3}}
 \uput[r](3,2){\scriptsize {2}}
 \uput[r](3,0){\scriptsize {4}}
 \psline{<-}(1.1,1.1)(2.9,2)
 \psline{<-}(1.1,1)(2.9,1)
 \psline{<-}(1.1,0.9)(2.9,0.1)
   \end{pspicture}
\end{center}
Then under the deg-lex order defined in \cite{Abdu1} and
\cite{Abdu2}, we get

\begin{theorem}(\cite{Abdu1,Abdu2})\
The set $S^{+c}$\ is a Gr\"{o}bner-Shirshov basis of the algebra
$U_q^+(G_2)$ (resp. $U^+(D_4$)).
\end{theorem}

   If we replace all $x$'s in $U_q^+(G_2)$ (resp. $U^+(D_4$))  with $y$'s, then we get a
   similar result for the negative part of the quantum group.
Then, by Theorem 2.7 in \cite{BoMa}, we have

\begin{theorem}(\cite{Abdu1,Abdu2})\ The set $S^{+c}\cup K\cup T \cup S^{-c}$ is a
    Gr\"{o}bner-Shirshov basis of the quantum group $U_q(G_2)$ (resp. $U^+(D_4)$).
\end{theorem}

\subsection{An embedding of recursively presented Lie algebras}

In 1961 G. Higman \cite{Higman} proved an important Embedding
Theorem which states that every recursively presented group can be
embedded in a finitely presented group. Recall that a group (or an
algebra) is called recursively presented if  it can be given by a
finite set of generators and a recursively enumerable set of
defining relations. If a group (or an algebra) can be given by
finite sets of generators and defining relations it is called
finitely presented. As a corollary to this theorem G. Higman proved
the existence of a~universal finitely presented group containing
every finitely presented group as a subgroup. In fact, its finitely
generated subgroups are exactly the finitely generated recursively
presented groups.

In \cite{Belyaev} V.Ya.Belyaev proved an analog of Higman's theorem
for associative algebras over a field which is a~finite extension of
its simple subfield. The proof was based on his theorem stating that
every recursively presented associative algebra over a field as
above can be embedded in a recursively presented associative algebra
with defining relations which are equalities of words of generators
and  $\alpha+\beta=\gamma$, where $\alpha, \beta, \gamma$ are
generators.

In recent paper \cite{Olshanskii} Y. Bahturin and A. Olshanskii, see
\S3.14, showed that such embedding can be performed distortion-free.
The idea of transition from  algebras to semigroups was also used by
G. P. Kukin in \cite {Kukin} (see also \cite{BokutKukin}).

The paper \cite{Chibrikov} appears as a byproduct of the author's
joint attempts with Y. Bahturin to prove Lie algebra analog of
Higman's Theorem. In particular, Y. Bahturin suggested to prove that
any recursively presented Lie algebra can be embedded in a Lie
algebra given by Lie relations of the type mentioned above. In the
paper \cite{Chibrikov}, the author shows that this is indeed true.
Namely, every recursively presented Lie algebra over a field which
is a~finite extension of its simple subfield can be embedded in a
recursively presented Lie algebra defined by relations which are
equalities of (nonassociative) words of generators and
$\alpha+\beta=\gamma$ ($\alpha, \beta, \gamma$ are generators). Note
that an existence of Higman's embedding for Lie algebras is still an
open problem (see \cite{Sapir}). It is worth to mention a result by
L. A. Bokut \cite{Bokut1} that for every recursively enumerable set
$M$ of positive integers, the Lie algebra
$$ L_M=\text{Lie}(a,b,c \, | \, [ab^nc]=0, n \in M),$$
where $[ab^nc]=[[\ldots [[ab]b] \ldots b]c]$, can be embedded into a
finitely presented Lie algebra.

By using Grobner-Shirshov basis theory for Lie algebras, in
\cite{Chibrikov}, the following theorem is given.

\begin{theorem}(\cite{Chibrikov})\
A recursively presented Lie algebra over a field which is a finite
extension of its simple subfield can be embedded into a recursively
presented Lie algebra defined by relations which are equalities of
(nonassociative) words of generators and  $\alpha+\beta=\gamma$,
where $\alpha, \beta, \gamma$ are generators.
\end{theorem}

\subsection{GSB for some monoids}

In the paper \cite{Turkey10},   Gr\"{o}bner-Shirshov bases for the
graph product, the Sch$\ddot{u}$tzenberger product and the
Bruck-Reilly extention of monoids are given, respectively, by using
Shirshov algorithm for associative algebras.

\end{document}